\documentclass[11pt]{amsart}

\usepackage{amssymb}
\usepackage{amsmath}
\usepackage{amscd}
\usepackage{graphicx}

\numberwithin{equation}{section}

\newtheorem{thm}{Theorem}
\newtheorem{lem}{Lemma}

\newtheorem{defn}{Definition}

\newtheorem{rem}{Remark}

\begin{document}

\title[Algebraic Markov equivalence for links in $3$--manifolds]{Algebraic Markov equivalence for links
in $3$--manifolds}

\author{S. Lambropoulou}
\address{ Department of Mathematics,
National Technical University of Athens,
Zografou Campus, GR-157 80 Athens, Greece.}
\email{sofia@math.ntua.gr}
\urladdr{http://www.math.ntua.gr/$\tilde{~}$sofia}

\author{C. P. Rourke}
\address{ Department of Mathematics, 
  University of Warwick, Coventry CV4 7AL, UK.}
\email{cpr@maths.warwick.ac.uk}
\urladdr{http://www.maths.warwick.ac.uk/$\tilde{~}$cpr}

\thanks{Most of this work was done while the first author held a position at the Laboratoire de
Math\'ematiques Nicolas Oresme, Universit\'e de Caen.}



\date{}
\maketitle

\begin{abstract}

Let $B_n$ denote the classical braid group on $n$ strands and let the {\em mixed braid group} $B_{m,n}$
be the subgroup of $B_{m+n}$ comprising braids for which the first $m$ strands form the identity
braid.  Let $B_{m,\infty}=\cup_nB_{m,n}$.  We will describe explicit algebraic moves on $B_{m,\infty}$
such that equivalence classes under these moves classify oriented links up to isotopy in a  link
complement or in a closed, connected, oriented $3$--manifold. The moves depend on a fixed link
representing the  manifold in  $S^3.$  More precisely, for link complements the moves are: the  two
familiar moves of the classical Markov equivalence 
together with {\em `twisted' conjugation} by certain 
loops $a_i$.  This means premultiplication by ${a_i}^{-1}$ and postmultiplication by a `combed' version
of $a_i$.  For closed $3$--manifolds there is an additional set of {\it `combed' band moves} which
correspond to sliding moves over the surgery link.  The main tool in the proofs is the one-move Markov
Theorem using  {\it $L$--moves} \cite{LR} (adding in-box crossings). The resulting algebraic
classification is a direct extension of the classical Markov Theorem that classifies links in $S^3$ up
to isotopy,  and potentially leads to powerful new link invariants, which have been explored in special
cases by the first author. It also provides a controlled range of isotopy moves,  useful for studying
skein modules of $3$--manifolds. 

\end{abstract}

\section{ Introduction and Overview }

By a classic result of H. Brunn and J.W. Alexander \cite{A}, \cite{Br}  any oriented knot in $S^3$ is
isotopic to the closure of a braid, and, by a theorem of A.A. Markov (and an improvement due to N.
Weinberg) \cite{Ma}, \cite{W}, \cite{Bi} there is a bijection (induced by `closing' the braid) between
isotopy classes of oriented links and equivalence classes of braids, the equivalence being generated by
braid isotopy and by two moves between braids:  {\it Markov conjugation} (conjugating by a
crossing) and the {\it Markov move} or {\it $M$--move} (adding an extra crossing at a rightmost
point).   In \cite{LR} we introduced a new type of braid move, the {\it $L$--move} (adding an in-box
crossing; view Figure 5 for abstract illustrations), and we showed that the equivalence relation
generated by $L$--moves and braid isotopy gives the same bijection.  Consequently, Markov conjugation and
Markov moves can be produced by $L$--moves (see Figures 12, 13).

\smallbreak 

The Markov Theorem can be regarded as a {\it geometrical\/} result (by thinking of braids as 
geometrical objects) or as an {\it algebraic\/} result (by thinking of braids as elements of the
classical braid group $B_n$). In the latter case the two moves  of the Markov equivalence have the
two well-known algebraic formulations. Similarly, the $L$--moves  have analogous algebraic formulations
(cf. \cite{LR}; Remark 2.2).

\bigbreak

Let now $V$ be the complement of a link. By `link complements'  we mean complements of both knots and
links and by `links' we always mean knots and links.  All links are considered oriented and piecewise
linear (PL), but will be mostly illustrated smooth for convenience.  By the Alexander Theorem, this link
is isotopic to the closure
$\widehat{B}$ of a braid~$B$. So, we can write  $V= S^3 \backslash \widehat{B}$ and $V$ can be
represented in $S^3$ by $\widehat{B}$. Further, let $V$ be  a closed, connected, oriented $3$--manifold
(we shall simply write ``closed $3$--manifold"). By classic results of Lickorish and Wallace \cite{Li},
\cite{Wa}
$V$  can be obtained from $S^3$ by surgery along a framed link with integral framings. Without loss of
generality the surgery link can be assumed to be the closure $\widehat{B}$ of a {\it surgery braid} $B$.
(Note that the framing of $\widehat{B}$ induces a framing on the surgery braid $B$.)
 So, we can write $V= \chi (S^3, \widehat B)$ and $V$ can be represented in $S^3$ by $\widehat{B}$. 
 Moreover, by the proof in \cite{Li}, all components of the surgery link can be assumed unknotted and, 
 as can be easily seen, they can be isotoped to the closure of a pure braid.
Thus, {\it for closed $3$--manifolds we may assume $B$ to be a pure braid\/.}  

\bigbreak

Let now $L$ be an oriented link in $V= S^3 \backslash \widehat{B}$ or $\chi (S^3, \widehat B)$. 
 Fixing $\widehat{B}$ pointwise we may represent $L$  in $S^3$ unambiguously by the {\it mixed link\/}
$\widehat{B}\bigcup L$, which consists of the fixed part $\widehat{B}$ and the `moving' part $L$ that
links with $\widehat{B}$ (see Figure 1 for an example). A {\it mixed link diagram\/} is a diagram
$\widehat{B} \bigcup \widetilde{L}$ of $\widehat{B}\bigcup L$ on the plane of $\widehat{B}$. This plane
is equipped with the top-to-bottom direction of $B$. By the Alexander Theorem and as explained in
\cite{LR} (cf. Theorem 5.3), a diagram $\widehat B \bigcup
\widetilde{L}$ of $\widehat{B}\bigcup L$ may  be turned  into a {\it mixed braid\/} $B \bigcup\beta$ 
with isotopic closure. (The closure of a braid is obtained by joining each pair of corresponding
endpoints by a simple arc.) This is a braid in $S^3$ with two different sets of strands, abstractly 
represented by a braid box with two differently coloured sets of strands.  The point here is that one of
the two sets comprises the {\it fixed subbraid $B$\/} and not any other Markov equivalent one. The
other set of strands representing the link in the manifold  $V$ is called the {\it  moving subbraid\/}.
See  Figure 1 for an example.  So, $V$ may be represented  in $S^3$ by the open braid $B$.

\bigbreak 

\begin{figure}[h]
\begin{center}
\includegraphics[width=4.5in]{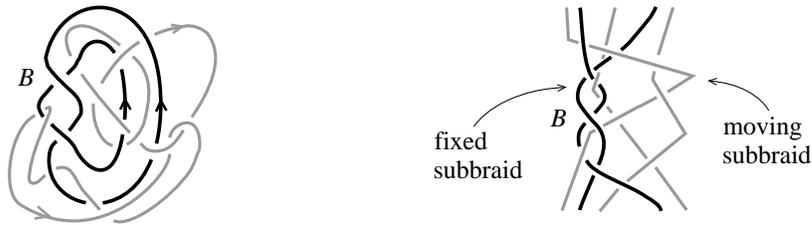}
\end{center}
\caption{ A mixed link and a mixed braid }
\label{figure1}
\end{figure}




Consider now an isotopy of $L$ in $V$. It follows from standard  results  of
PL Topology that $L_1$ and $L_2$ are two instances  of an isotopy in $S^3 \backslash \widehat{B}$ if and
only if the corresponding  mixed links $\widehat{B}\bigcup L_1$ and $\widehat{B}\bigcup L_2$ are
isotopic in $S^3$ by an ambient isotopy which keeps $\widehat{B}$ pointwise fixed. See
\cite{RS}. In terms of diagrams, the mixed link isotopy will not involve Reidemeister moves of the fixed
part. 

\smallbreak

The first stage of surgery along a framed link $\widehat{B}$ is to pass from $S^3$ to the link complement
$S^3 \backslash \widehat{B}$. Thus, an isotopy of $L$ in $\chi (S^3, \widehat B)$ can be viewed as an
isotopy in $S^3 \backslash \widehat{B}$, but  with the extra freedom for  $L$ to slide  across the disc
that the parallel curve of a framed component of $\widehat B$ bounds in $\chi (S^3, \widehat B)$. This
isotopy move is similar to the second move of the Kirby calculus. As noted in \cite{LR}, the first part
of the move is just isotopy in $S^3 \backslash \widehat B$, so we only need to consider the essential
part, where a little band of $L$ very close to the surgery component slides over the component,
according to the framing and orientation conventions. We shall call this move a {\it band move}.   A
band move takes place in an arbitrarily thin tubular neighbourhood of the component of the surgery link,
so by `band move' we may unambiguously refer to both the move in the $3$--space and its projection on
the plane of $ \widehat B$.  In terms of diagrams, the mixed link equivalence in $S^3$ includes the band
moves (two types, depending on the orientation of the little band, which are related by a twist of
the little band; see Figure~2). For more details the reader is referred to \cite{LR}, Theorems 5.2 and
5.8. 

\bigbreak 

\begin{figure}[h]
\begin{center}
\includegraphics[width=5.5in]{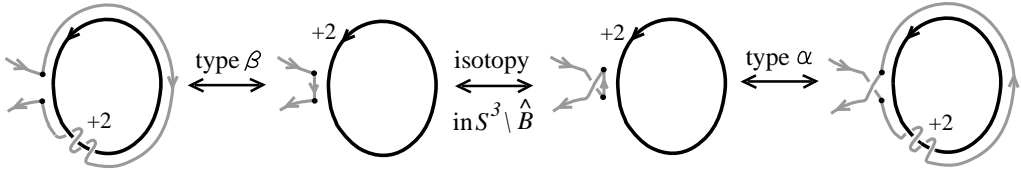}
\end{center}
\caption{ The two types of band moves and their relation }
\label{figure2}
\end{figure}

Let's see now how the mixed link isotopy translates on the level of mixed braids.

\begin{defn}{\rm \  A  {\it braid band move\/} is a move between mixed braids, which is a band
move between their closures. It starts with a little band oriented downwards, which, before sliding along
a surgery strand, gets one twist {\it positive\/} or {\it negative\/}.  View  Figure 3. In the sequel
we shall omit the word `braid' and we shall just say `band move'.  }
\end{defn}

\bigbreak 

\begin{figure}[h]
\begin{center}
\includegraphics[width=3.5in]{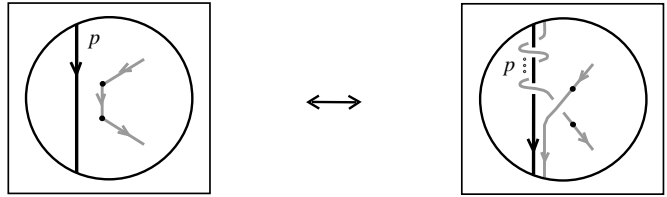}
\end{center}
\caption{ The  band move for mixed braids }
\label{figure3}
\end{figure}

\begin{defn}[$L$--moves for mixed braids]{\rm \ Let $B \bigcup \beta$ be a mixed braid in $S^3$  and $P$
a point of an arc of the {\it moving subbraid\/} $\beta,$ such that $P$ is not vertically aligned with
any  crossing. Doing an {\it $L$--move\/}  at $P$ means breaking the arc at $P$, bending the two
resulting smaller  arcs slightly apart  by a  small isotopy and stretching them vertically, the upper
downwards and the lower upwards, and both {\it over\/}  or {\it under\/} all other arcs of the diagram,
so as to introduce two new corresponding strands with endpoints on the vertical line of $P$. Stretching
the new strands over  will give rise to an {\it $L_o$--move\/} and under to an {\it  $L_u$--move\/}. 
See Figure 4. }
\end{defn}

\bigbreak 

\begin{figure}[h]
\begin{center}
\includegraphics[width=5.3in]{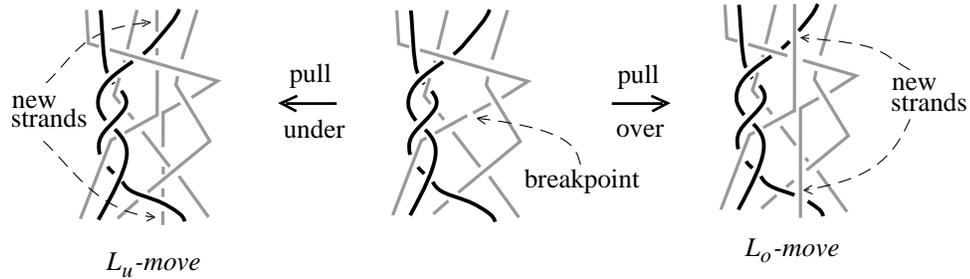}
\end{center}
\caption{  The two types of $L$--moves }
\label{figure4}
\end{figure}

\noindent Using a small braid isotopy, an $L$--move can be equivalently seen with a crossing (positive or
negative) formed. See Figure 5. 

\bigbreak 

\begin{figure}[h]
\begin{center}
\includegraphics[width=4.5in]{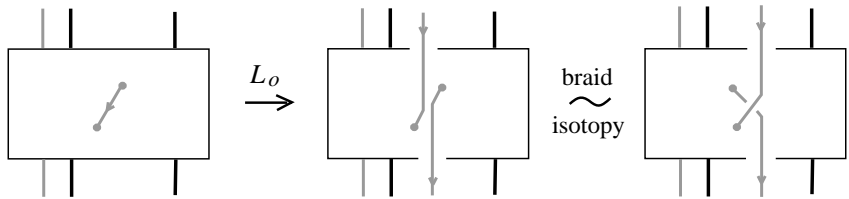}
\end{center}
\caption{  An $L$--move introduces a crossing }
\label{figure5}
\end{figure}

\noindent Clearly, two mixed braids that differ by
an  $L$--move have isotopic closures, since the $L$--move corresponds to introducing a twist
in the mixed link. $L$--moves and mixed braid isotopy generate an equivalence  relation on
mixed braids called  {\it $L$--equivalence}. 
Our method of  proving the one--move (and  the classical) Markov Theorem (Theorem 2.3 in
\cite{LR})  ensures that the arcs of the diagram that are oriented downwards do not
participate in the proof. This led us to the following result (Theorem 5.5 and Theorem 5.10), which is
our starting point in this paper.  

\begin{thm}[Geometric Markov Theorem for $V= S^3 \backslash \widehat B$ or $\chi (S^3, \widehat B)$] \
Two  oriented links in $S^3 \backslash \widehat B$ are isotopic if and only if  any two corresponding
mixed braids in $S^3$ differ by mixed braid isotopy and a finite sequence of $L$--moves that do not
touch the fixed subbraid $B$.  

Moreover, if the two links  lie in $\chi (S^3, \widehat B)$, the mixed braids differ by mixed braid
isotopy, by $L$--moves that do not touch the fixed subbraid $B$ and by braid band moves.
\end{thm}

 The paper  is concerned with the corresponding algebraic formulation. 

\bigbreak 

\begin{figure}[h]
\begin{center}
\includegraphics[width=1.4in]{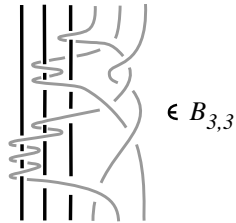}
\end{center}
\caption{ An algebraic mixed braid }
\label{figure6}
\end{figure}

The  braid structures for links in these  manifolds (as well as for links in handlebodies) have been
established and studied in \cite{L2}. These are either the extended  braid groups $B_{m,n}$, whose
elements are called {\it algebraic mixed braids} and they have the first
$m$  strands  forming the identity subbraid (view Figure 6 for an example), or appropriate cosets
$C_{m,n}$ of these groups, that depend on the specific manifold. 
  More precisely, $B_{m,n}$ has the presentation:
\[ 
B_{m,n} = \left< \begin{array}{ll}  \begin{array}{l} 
a_1, \ldots, a_m,  \\ 
\sigma_1, \ldots ,\sigma_{n-1}  \\
\end{array} &
\left|
\begin{array}{l} \sigma_k \sigma_j=\sigma_j \sigma_k, \ \ |k-j|>1   \\ 
\sigma_k \sigma_{k+1} \sigma_k = \sigma_{k+1} \sigma_k \sigma_{k+1}, \ \  1 \leq k \leq n-1  \\
{a_i} \sigma_k = \sigma_k {a_i}, \ \ k \geq 2, \   1 \leq i \leq m,    \\ 
 {a_i} \sigma_1 {a_i} \sigma_1 = \sigma_1 {a_i} \sigma_1 {a_i}, \ \ 1 \leq i \leq m  \\
 {a_i} (\sigma_1 {a_r} {\sigma_1}^{-1}) =  (\sigma_1 {a_r} {\sigma_1}^{-1})  {a_i}, \ \ r < i.  
\end{array} \right.  \end{array} \right>  
\]
where the generators $a_i$ and $\sigma_j$ are as illustrated in Figure 7.

\bigbreak 

\begin{figure}[h]
\begin{center}
\includegraphics[width=5in]{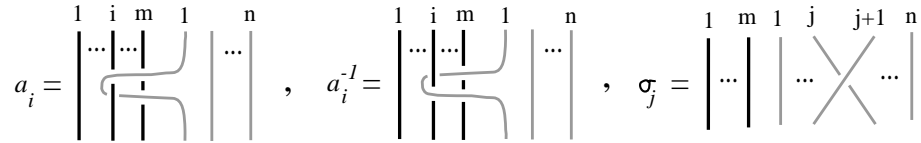}
\end{center}
\caption{ The `loops' $a_i, \ {a_i}^{-1}$ and the crossings $\sigma_j$ }
\label{figure7}
\end{figure}

The groups $B_{m,n}$ are the appropriate braid structures for studying knots and links in the complement
of the $m$--unlink or a connected sum of $m$ lens spaces of type $L(p,1)$ or a handlebody of genus
$m$. For the first two cases of manifolds it is easy to formulate the analogue
 of Markov Theorem algebraically (see first two examples of Section 4). In \cite{HL} an algebraic
formulation of the Markov Theorem for handlebodies was proven in terms of the groups $B_{m,n}$, one
version using algebraic
$L$--moves and another using Markov equivalence (cf. \cite{HL}, Theorems 4 and 5).  In that case there
was no surgery involved. The conceptual difficulty there was related  to the fact that conjugations by
the
$a_i$'s were not permitted. 

\bigbreak 

For a generic $V= S^3\backslash \widehat B$ or $\chi (S^3, \widehat B)$ the fixed subbraid $B$ is not the
identity braid. {\it Parting} a mixed braid means to separate its endpoints into two different sets, so
that the resulting braids  have isotopic closures. Figure 9 illustrates different partings of an
abstract mixed braid. {\it Combing} a parted mixed braid means to separate the fixed subbraid from the
moving part, using mixed braid isotopy. See Figure 8 for an abstract illustration. These operations are
discussed in detail in Section 2 and Section 3. By parting  and  combing mixed braids, it was shown in
\cite{L2}, Section 6 that knots and links in $V$ may be  represented by mixed braids in the groups
$B_{m,n}$ followed by the natural embeddings of $B$ in the groups $B_{m+n}$.  Hence, that the
braid structures related to $V$ are the cosets of the subgroups
$B_{m,n}$  in the groups $B_{m+n}$ ($n \in {\Bbb N}$), containing the embedded fixed subbraid $B$.
                                                                        
\bigbreak 

\begin{figure}[h]
\begin{center}
\includegraphics[width=4.4in]{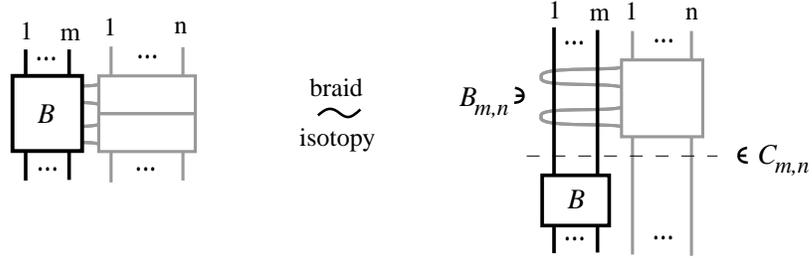}
\end{center}
\caption{ An abstract parted and combed mixed braid }
\label{figure8}
\end{figure}

\noindent {\bf The main results:} The main results of this paper are Theorem 4 and Theorem 5 in Section
3. Theorem 4 gives the algebraic braid equivalence of combed mixed braids for knot
complements  and Theorem 5 gives the algebraic braid equivalence of combed mixed braids for 
closed $3$--manifolds. Our strategy for proving these Theorems is
the following. We first part the mixed braids and we translate the  $L$--equivalence and the braid band
moves of Theorem~1 to an {\it equivalence of parted mixed braids}.  Here the generators of the groups
$B_{m,n}$ become apparent in the equivalence. Also, the braid band moves assume a special form.  
 These are done in Section 2. See Theorem 2 for link complements,  and Lemma 5 and Theorem 3 for closed
$3$--manifolds.

\smallbreak
In Section 3 we comb the parted mixed braids and we translate the parted mixed braid equivalence to an
{\it equivalence of algebraic mixed braids}. For both link complements and closed $3$--manifolds, 
Markov move and Markov conjugation remain equivalence moves between combed and algebraic  mixed braids.
But in place of the conjugation by a `loop' $a_i$ we need to introduce the {\it twisted
conjugation}, which takes into account the combing of the loop through the fixed subbraid: 
\[
\beta  \sim  {a_i}^{\mp 1} \beta {\rho_i}^{\pm 1} 
\]
where $\rho_i$ is the combing of the loop  $a_i$ through $B,$ for \ $\beta, a_i, \rho_i \in B_{m,n}$.
 See Figure 19 for illustrations. 

\smallbreak
Moreover, a {\it parted band move} after combing is the {\it composition of an algebraic band  move with
the combing of the parallel strand} through the surgery subbraid. An {\it algebraic band move} is a braid
band move between elements of the groups $B_{m,n}$ and it has the algebraic expression: 
\[ 
\beta_1 \beta_2 \ \sim \ \beta'_1 \, {t_{k,n}}^{p_k} \,
{\sigma_n}^{\pm 1}\, \beta'_2
\]   
\noindent where $\beta_1, \beta_2 \in B_{m,n}, \ t_{k,n}$ is a Markov conjugate of the loop $a_k, \ p_k
\in {\Bbb Z}$ is the framing of the $k$th surgery component of the surgery link and  $\beta'_1, \beta'_2 
\in B_{m,n+1}$  are the words $\beta_1, \beta_2$, but with certain substitutions that indicate the
pulling of the parallel strand to the right of the braid. See Definition 7 and Figure 20. Then, a
{\it combed band move} has the algebraic expression:
\[
 \beta_1 \beta_2 \ \sim \ \beta'_1 \, {t_{k,n}}^{p_k} \,
{\sigma_n}^{\pm 1}\, \beta'_2 \, r_k
\]
\noindent where $r_k$  is the  combing  of the parted parallel strand  to the $k$th surgery strand
through the surgery braid. For an illustration here see Figure 21. 

\smallbreak

Finally, in Section 4 we give explicit examples, including complements of daisy chains, the lens
spaces, the complement of the Borromean rings or a closed $3$--manifold obtained by surgery along them.
We also discuss the case where the surgery braid is not a pure braid (Lemma 9) and we present as an
example the complement of a trefoil or a manifold obtained by surgery along it.

\bigbreak

 This paper is sequel to \cite{LR}, \cite{L2} and \cite{HL}. It sets out the necessary algebraic
formalism for constructing knot invariants in $3$--manifolds using braid machinery,
for example  via constructing Markov traces on appropriate  algebras, quotients of the group algebras of
the  braid groups  $B_{m,n}$. (See  \cite{J} for the classical case of links in $S^3$). In the case
$m=1$, $B_{1,n}$ is the Artin group of type $\mathcal B$. See, for example, \cite{L1} and references
therein for the construction of the analogues of the 2--variable Jones polynomial (homflypt) for links in
the solid torus. The case of $L(p,1)$ is being studied in \cite{LP}. Theorem 5 gives a very good
control over the  band moves of links in closed $3$--manifolds, and this is very useful for the study of 
skein modules of $3$--manifolds. For skein modules of $3$--manifolds see, for example, \cite{P} and
references therein. 

\smallbreak

A final comment is now due. In our set--up  the  manifold is represented in  $S^3$ by  a fixed
link. In the case of a closed $3$--manifold the surgery link is not unique up to isotopy. In fact, it may
be altered via the Kirby calculus \cite{Ki}. Then the corresponding surgery braids are equivalent under
moves described in \cite{KS}. One could then consider combining the mixed braid equivalence given in the
present paper with the braid equivalence of \cite{KS}. This will be the subject of a subsequent paper.

\section{Markov equivalence for parted mixed braids}

This section is an intermediate step towards the algebrization of Theorem~1. Here the
mixed  braids resulting from the mixed links have all, say $m$, strands of the fixed
part $B$  occupying the first $m$ positions of the mixed braid. 

\begin{defn}{\rm \ A {\it parted mixed braid\/} is a mixed braid $B \bigcup\beta$ on $m+n$
strands, such that the first $m$ endpoints are those of the subbraid $B$ and the last $n$
endpoints are those of $\beta.$  Parted mixed braids are denoted in the same way as mixed
braids. We number their fixed strands from $1$ up to $m$ and their moving
strands from $1$ up to $n.$  (See left hand side of Figure 14 for an abstract illustration.)   } 
\end{defn}

\begin{lem}{ \  Every mixed braid may be represented by a parted mixed braid with
isotopic closure (cf. Section 6 in \cite{L2}, compare with Lemma 1 in \cite{HL}).
 }
\end{lem}

\begin{proof}  Indeed, let  $B \bigcup\beta$ be a
mixed braid. To see this we simply attach arbitrarily arrays of labels `o' or `u' to 
corresponding pairs of endpoints of the moving subbraid $\beta$, with as many entries as the number of
fixed strands on their right,  and  we pull the strands of  corresponding endpoints to the right, {\it
over\/} or {\it under\/} each  strand of $B$ that lies on their right, according to the  label in the
array of the pulled strands. We start from the rightmost pair respecting   the position of the
endpoints. View the  first two illustrations of Figure 9 for the parting of an abstract mixed
braid. Obviously, the closures of the initial and of the parted mixed braid are isotopic (they differ by
planar isotopy and by mixed Reidemeister II moves). 
\end{proof} 

\bigbreak 

\begin{figure}[h]
\begin{center}
\includegraphics[width=5.5in]{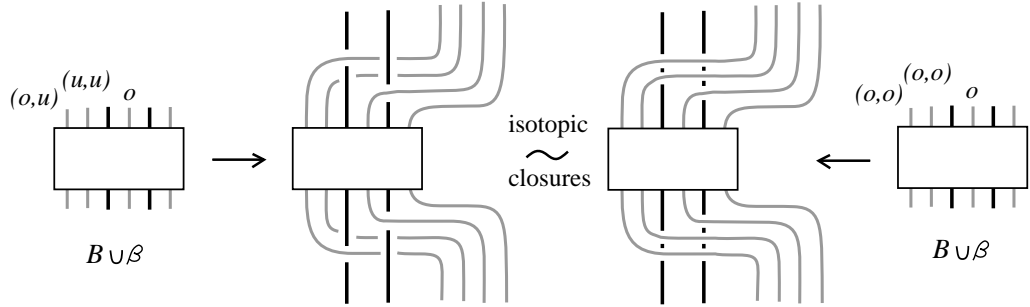}
\end{center}
\caption{ Parting a mixed braid -- the standard parting }
\label{figure9}
\end{figure}

It follows that two different partings of a mixed braid give rise --upon closure-- to isotopic 
 mixed links. For this reason we could fix the pulling of the moving
strands during the parting process to be always  over or always under each
strand of $B$ that lies on its right. If the pulling is always {\it over\/} we shall refer
to it as {\it the standard parting.\/}  See the two right illustrations of  Figure 9.  

\smallbreak 

Pulling a moving strand `under' a fixed strand instead of `over', it simply corresponds to the fact that
the closure of the moving strand  crosses a hypothetical closing arc $k$ of the fixed subbraid $B$, 
and this is an allowed isotopy move in the manifolds considered here. See Figure 10. (We note
 that this is not true in the case of handlebodies and this is the reason why conjugation by the $a_i$'s
is not permitted; see \cite{HL} for a detailed analysis.)   

Lemma 2  below gives the  relation of an
arbitrary parting with the corresponding standard parting, and it is very instructive, as it brings the 
`loops'  $a_i$ into the parted braid equivalence. Note that the  elementary algebraic mixed braids $a_i$
and their inverses, together with the crossings $\sigma_j$ (all defined in  Figure~7) are
clearly the {\it geometric generators} of the moving part of a parted mixed braid. 

\bigbreak 

\begin{figure}[h]
\begin{center}
\includegraphics[width=5.4in]{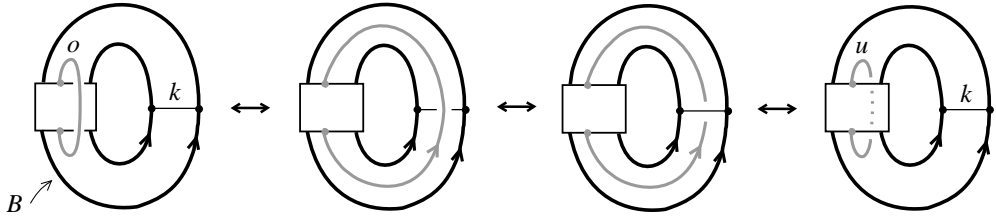}
\end{center}
\caption{ Crossing the arc $k$ }
\label{figure10}
\end{figure}

 Let now $C_{m,n}$ denote the set of
parted mixed braids  on  $n$ moving strands related to $V= S^3 \backslash \widehat B$ or $\chi (S^3,
\widehat B)$. By adding an extra moving strand on the right of a parted mixed braid on  $n$ moving
strands, $C_{m,n}$ embeds naturally into $C_{m,n+1}$.  Let 
$C_{m,\infty} := \bigcup_{n=1}^{\infty} C_{m,n}$ denote the disjoint union of  all sets
$C_{m,n}$. We define below some moves in  $C_{m,\infty}$.

\begin{defn}{\rm 
{\it (1) \ Loop conjugation\/} of a parted mixed braid in $C_{m,n}$ is
its  concatenation from above by a loop $a_i$  (or by ${a_i}^{-1}$) and from
below by  ${a_i}^{-1}$ (corr.  $a_i$). 

\vspace{.05in}
\noindent {\it (2) \  Markov conjugation\/} of a parted mixed braid  in $C_{m,n}$ is 
its concatenation from above  by a crossings $\sigma_j$ (or by ${\sigma_j}^{-1}$) and from below by
${\sigma_j}^{-1}$  (corr.  $\sigma_j$).

\vspace{.05in}
\noindent {\it (3)} \ A {\it parted $L$--move\/} is defined to be  an  $L$--move  between
parted mixed braids. (See left hand illustration of Figure 15.)  

\vspace{.05in}
\noindent  {\it (4)} \ An  {\it $M$--move} is the insertion of a crossing
${\sigma_n}^{\pm 1}$  at the right hand side of a parted mixed braid on $n$
moving strands. Undoing an $M$--move is the reverse operation. See Figure 14. 
 }
\end{defn}

\begin{lem}{ \ Consider a mixed braid on $m+n$ strands and an arbitrary parting of it in
 $C_{m,n}.$ Then, up to Markov conjugation, this parting differs
from its corresponding standard parting by a finite sequence of loop conjugations. }
\end{lem}

\begin{proof}
By an inductive argument we may assume that all moving strands
 from the $1$st up to the $(j-1)$st  are pulled `over' all the fixed strands that lie on
their right. Consider now the $j$th moving strand. Upon parting, this lands on the $j$th
position of the moving part of the resulting parted mixed braid. See Figure 11. By a mixed
braid isotopy  we bring the $j$th moving strand on top of the other moving strands.  Note
that the braid isotopy is independent of the parting labels attached to the $j$th moving
strand. Then by Markov conjugation by the word $(\sigma_1\ldots \sigma_{j-1})$ we bring the
$j$th moving strand to the first position of the moving subbraid. See top row of Figure 11.

 \bigbreak 

\begin{figure}[h]
\begin{center}
\includegraphics[width=5.7in]{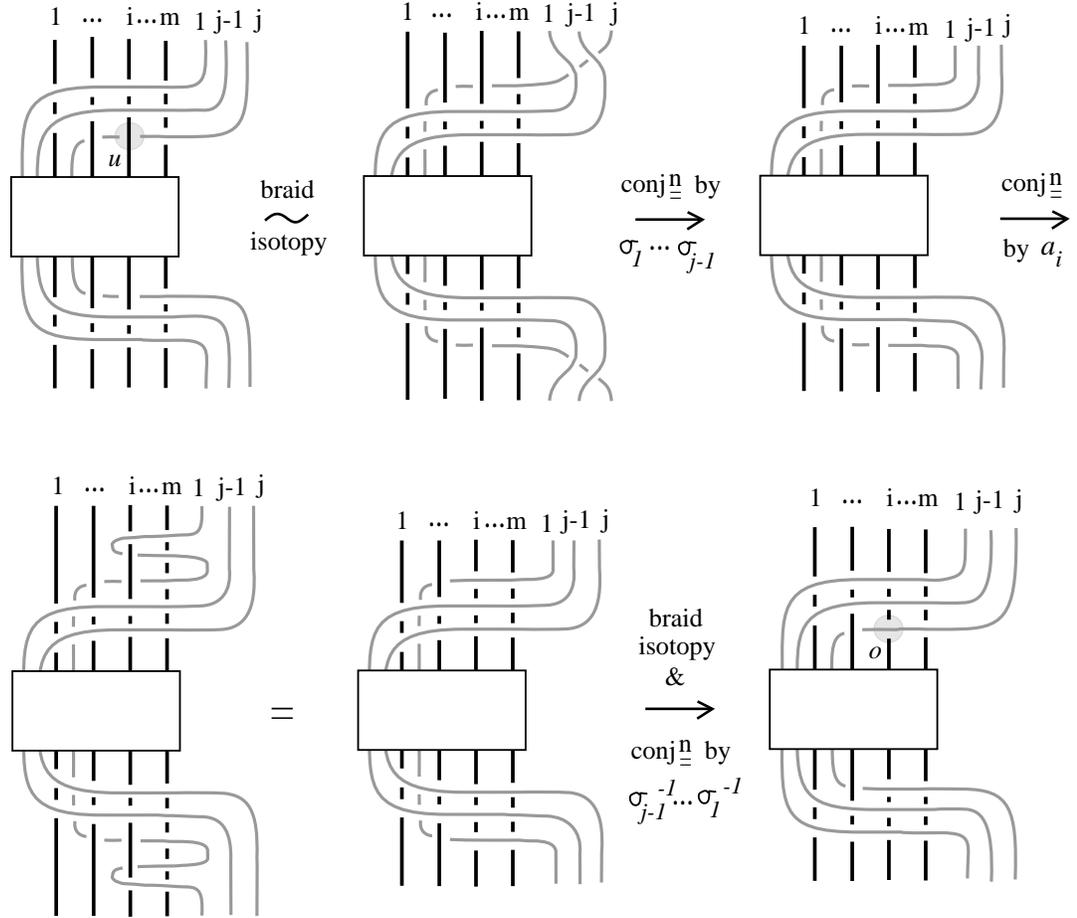}
\end{center}
\caption{ Change of parting labels $\longleftrightarrow$ conjugation by $a_i$ }
\label{figure11}
\end{figure}

\noindent  Let now the parting label of the original $j$th moving strand for the $i$th
fixed strand be  `under'. By a second inductive argument we may  assume  that the parting
labels of the original $j$th  strand are all `over' for the $(i+1)$st up to the $m$th
fixed strand. Then, 
conjugation by $a_i$ changes the label `under' to `over'. By applying once more mixed braid
isotopy and Markov conjugation by the word $({\sigma_{j-1}}^{-1}\ldots 
{\sigma_1}^{-1})$  we obtain a parted mixed braid identical to the initial one except for
the place of the one crossing in question, which is switched.  See second row of
Figure 11. Continuing backwards with the remaining parting labels of the
$j$th moving strand we change them all  in this manner  to `over', and this ends the proof.
\end{proof}

\begin{rem}{\rm \ It follows from the proof of Lemma 2 that changing a parting label from `under' to
`over' corresponds in  $C_{m,\infty}$ to conjugation by some $a_i$. 
 }
\end{rem}

\begin{lem}{ \ A  mixed braid with  an $L$--move performed can be parted to a parted mixed
braid with a parted $L$--move performed. (Compare with Lemma 2 in \cite{HL}.) 
 }
\end{lem}

\begin{proof}
 If the $L$--move is an 
$L_o$--move we part its strands by pulling them  to the right and over all other strands
in between.  Then the crossing of the 
$L$--move slides over to the right by a braid isotopy. See Figure~12. The case of an
$L_u$--move is  analogous: here we pull the two strands  under the fixed 
strands in between.  
\end{proof}

\bigbreak 

\begin{figure}[h]
\begin{center}
\includegraphics[width=4.5in]{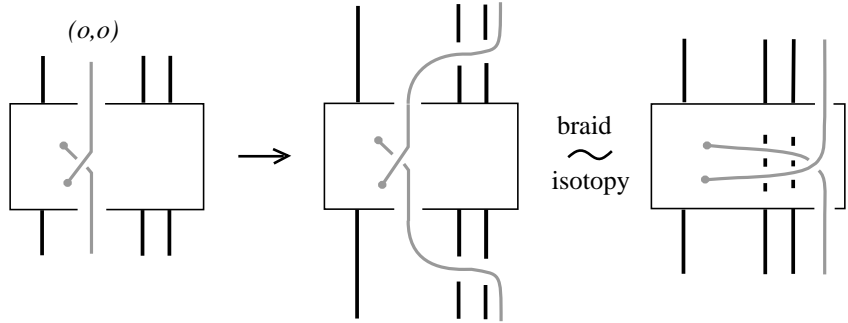}
\end{center}
\caption{  Sliding an $L_o$--move to the right }
\label{figure12}
\end{figure}

\begin{lem}{ \ Markov conjugation and  $M$--moves can be realized by a sequence of parted $L$--moves.
Conversely, a parted $L$--move is a composition of an $M$--move and Markov conjugation.
 }
\end{lem}

\begin{proof}
It is clear that an $M$--move is a special case of a parted $L$--move. The one-move Markov Theorem in
$S^3$  implies that Markov conjugation in $S^3$ can be realized by a sequence of $L$--moves (cf.\
4.1 in \cite{LR}). The same arguments apply to both link complements and closed $3$--manifolds. 
 But we would like to give a second direct proof of Lemma 4, which is an 
adaptation  for the case of knot complements and closed $3$--manifolds of a direct proof
for the classical case of $S^3$,  given by Reinhard H\"aring-Oldenburg.  In Figure
13 we start with a parted mixed braid conjugated by $\sigma_j$. After performing an
appropriate parted $L_o$--move, braid isotopy and  undoing  another parted $L_o$--move we end up
with the original  mixed braid free of conjugation by the $\sigma_j.$
 Conversely, as it becomes clear from Figure 15, a parted $L$--move is a composition of an $M$--move and
Markov conjugation.   
\end{proof}

\bigbreak 

\begin{figure}[h]
\begin{center}
\includegraphics[width=5.7in]{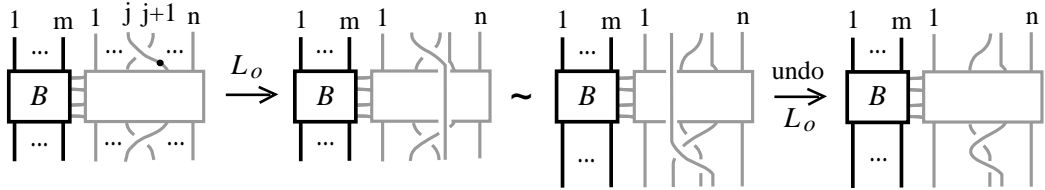}
\end{center}
\caption{  Conjugation by ${\sigma_j}^{-1}$ is a composition of $L$--moves }
\label{figure13}
\end{figure}

\bigbreak 

\begin{figure}[h]
\begin{center}
\includegraphics[width=3.5in]{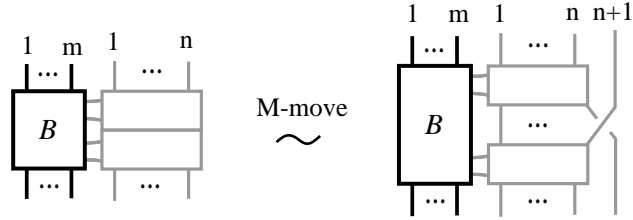}
\end{center}
\caption{ The $M$--move }
\label{figure14}
\end{figure}

\bigbreak 

\begin{figure}[h]
\begin{center}
\includegraphics[width=3.5in]{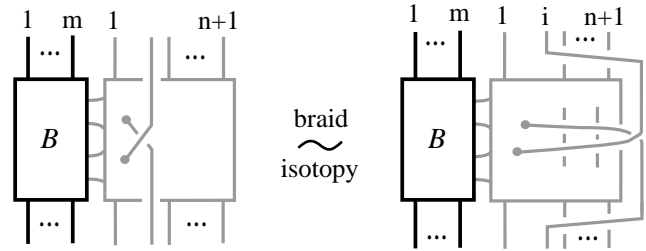}
\end{center}
\caption{ Parted $L$--move $\longleftrightarrow$ $M$--move and Markov conjugation }
\label{figure15}
\end{figure}

We are now in a position to state two versions of the analogue of the Markov Theorem for
parted mixed braids in $S^3 \backslash \widehat B$. 

\begin{thm}[Parted Version of Markov Theorem for $V = S^3 \backslash \widehat
B$] \ Two  oriented links in $S^3 \backslash \widehat B$ are isotopic if and only if 
any two corresponding parted mixed braids in $S^3$ differ by a finite sequence of parted
$L$--moves and loop conjugations.

\vspace{.03in}
Equivalently, two  oriented links in $S^3 \backslash \widehat B$ are isotopic if and only
if  any two corresponding parted mixed braids in $C_{m,\infty}$ differ by a  finite
sequence of  $M$~--~moves,   Markov conjugations and loop conjugations. 
\end{thm} 

\begin{proof}
It follows immediately from Theorem 1 for $V= S^3 \backslash \widehat B$ and from Lemmas 2, 3 and 4.
\end{proof}

We would like to extend Theorem 2 to parted mixed braids in closed $3$--manifolds. Lemma 5
below sharpens the band moves of Theorem 1 for $\chi (S^3, \widehat B)$ and it shows the
effect of  parting  on band moves. 

\begin{defn}{\rm \ {\it A parted band move\/} is defined to be a  band move
between parted mixed braids, such that: it takes place at the top part of the braid (before any surgery
crossings are encountered) and the little band {\it starts from the last strand} of the moving
subbraid and it {\it moves over\/}  each component of the surgery braid, until it reaches from the right
the specific component. After the band move is performed we apply to the resulting mixed braid
 the standard parting.
  } \end{defn}

\noindent  See Figure 16 for an example of a positive parted band move, where the moving
part has been simplified to the identity braid. 

\bigbreak 

\begin{figure}[h]
\begin{center}
\includegraphics[width=3.5in]{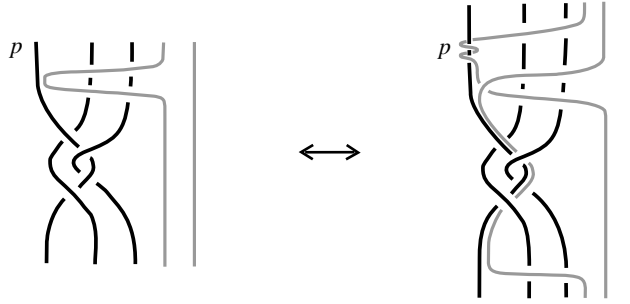}
\end{center}
\caption{  A  parted band move }
\label{figure16}
\end{figure}

\begin{lem}{ \ A band move may be always assumed, up to $L$--equivalence, to take
place at the top part of a mixed braid and on the right of the specific surgery strand.

Moreover, performing a band move on a parted mixed braid and then parting, the
result is equivalent, up to $L$--moves and loop conjugation, to performing a parted band
move.
 }
\end{lem}

\begin{proof}
  In a mixed braid $B \bigcup\beta_1$ consider a little band that
has approached a specific  surgery strand of $B$ from the right and is about to
perform a band move.   Pull the little band up to the top along the surgery strand.
See  illustrations 1 and 3 of Figure 17. 
 Then do braiding to obtain a mixed braid $B \bigcup\beta_2.$ This is  $L$--equivalent to
$B \bigcup\beta_1.$  Note that the edge arc of the little band is still
there in $B \bigcup\beta_2,$ because it is a down-arc. Now, using this arc, perform a {\it
top band move\/} in $B \bigcup\beta_2$ and  call the resulting mixed braid $B
\bigcup\alpha_2.$ See illustration 4 of Figure 17. Let also $B
\bigcup\alpha_1$ be the mixed braid obtained after performing the band move in the first
mixed braid
$B\bigcup\beta_1.$ See illustration 2 of Figure 17, but consider only the braid between
the dotted lines. Then 
$B \bigcup\alpha_1$ differs from $B \bigcup\alpha_2$ by exactly the same sequence of
$L$--moves, $L_1,\ldots,L_k$ say, that separate $B \bigcup\beta_1$ and $B
\bigcup\beta_2,$ since the isotopies separating the corresponding closures are identical. 
Compare the corresponding diagrams of Figure 17. Thus we showed that:

 \[ B \bigcup\beta_1  \stackrel{general \, band \, move}{\sim}  B \bigcup\alpha_1
\]  
 \[\Longleftrightarrow \]  
\[
 B \bigcup\beta_1  \stackrel{L_1,\ldots,L_k}{\sim}  B \bigcup\beta_2 
\stackrel{top \, band \, move}{\sim}  B
\bigcup\alpha_2  \stackrel{{L_1}^{-1},\ldots,{L_k}^{-1}}{\sim}  B \bigcup\alpha_1.  
\]  

\bigbreak 

\begin{figure}[h]
\begin{center}
\includegraphics[width=14cm]{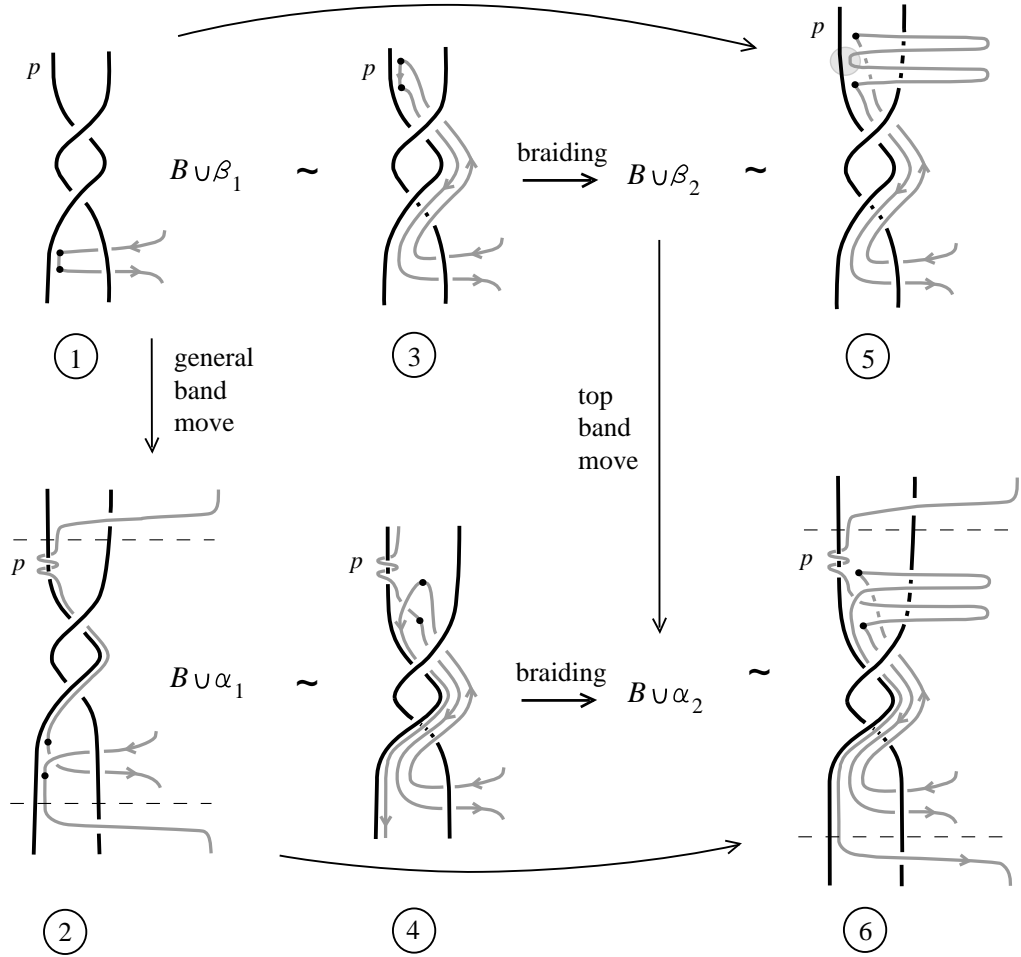}
\end{center}
\caption{  The proof of Lemma 5 }
\label{figure17}
\end{figure}

The first statement of the Lemma is proved. Consider now the same setting as above, but with $B
\bigcup\beta_1$ being a parted mixed braid. Perform on it a  band move and part the resulting new mixed
braid by pulling the two new strands over all strands in between to the last position of the moving
part. See  illustration 2 of Figure 17, where the parting is now included. At the same time pull the
little band of $B \bigcup\beta_1$ up to the top and then horizontally to the right, {\it over\/} all 
strands in between, until it reaches the last position of the moving part.  Then do a similar pull-back
to the left up to the specified surgery strand.  See 
illustration 5 of Figure 17. Now perform a parted band move at the place marked with a shaded disc. See 
illustration 6 of Figure 17. As above, this last mixed braid operation does not create any
new up-arcs and it does not interfere with the band move. Finally, part by the standard
parting the new strands created by braiding the up-arcs from the pulling along the surgery strand,  and
part last the new strands created from the parted band move, by pulling them over all strands in between
to the last position of the moving part. 
 Clearly, the two pairs (related to $(1,5)$ and $(2,6)$ of Figure 17) of parted
mixed braids involved differ by the same parted 
$L$--moves together with loop conjugation that comes from the parting. 

\smallbreak

  Further note that, by braid isotopy and loop conjugation, the $p$ twists of a general
band move may take place anywhere along the surgery strand,  so just as well at the
top, as in Definition 5 of a parted band move.  Thus, we showed that a band move on a parted mixed braid
can be accomplished with a parted band move, up to $L$--moves and loop conjugation.
\smallbreak

Finally, if the little band lies on the left of the surgery component we pull it
horizontally over the surgery strand and
to the right and then we pull it slightly back to the left, so that it approaches the
surgery strand from the right. Up to here it is only braid isotopy. We now perform a 
band move before and a band move after  and we notice  that the two final mixed braids differ
by conjugating a half twist of the framing. In any case, after parting the two final
mixed braids are the same.
\end{proof}

Then Theorem 1 for $\chi (S^3, \widehat B)$ and Lemma 5 extend Theorem 2 to the following.

\begin{thm}[Parted Version of Markov Theorem for $\chi (S^3, \widehat B)$] \ Two 
oriented links in  $\chi (S^3, \widehat B)$ are isotopic if and only if  any two
corresponding  parted mixed braids in $C_{m,\infty}$ differ by a finite sequence of parted
$L$--moves,  loop  conjugations and  parted band moves.

\vspace{.03in}
Equivalently, two oriented links in  $\chi (S^3, \widehat B)$ are isotopic if and only if 
any two corresponding  parted mixed braids in $C_{m,\infty}$ differ by a finite sequence
of $M$~--~moves,   Markov conjugation, loop conjugations and parted band moves.
\end{thm}

\section{Markov equivalence for combed and algebraic mixed braids}

In this section we explain the combing of parted mixed braids and we translate the equivalence of
Theorem 2 and Theorem 3 to equivalence of algebraic mixed braids. 

\bigbreak

Let  $V= S^3 \backslash \widehat B$ or $\chi (S^3, \widehat B).$  Unless $V$ is 
the complement of the $m$--unlink or a connected sum of $m$ lens spaces of  type
$L(p,1)$, where the fixed subbraid $B$  is the identity braid on $m$ strands, concatenating two
elements of $C_{m,n}$ is not a closed operation, since it alters the braid description of
the manifold.  So, the set $C_{m,n}$ of parted mixed braids is not a subgroup of  
$B_{m+n}$.  Yet, as shown in \cite{L2}, Section 6,  using Artin's combing for pure braids, the moving
part of a parted mixed braid can be combed away from the fixed subbraid, so that this
latter remains free of mixed linking at the bottom of the parted mixed braid. Thus, the parted mixed
braid splits into the concatenation of two parted mixed braids: the {\it `algebraic'
part\/} at the top, which has as fixed subbraid the identity braid on $m$ strands and the
{\it `coset'  part\/} at the bottom consisting of the fixed braid $B$ embedded naturally in  $B_{m+n}$.
See Figure 18 for an abstract illustration.  The result will be called a {\it
combed mixed braid}. 

\bigbreak 

\begin{figure}[h]
\begin{center}
\includegraphics[width=2in]{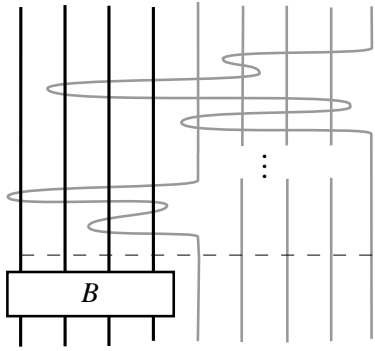}
\end{center}
\caption{ Artin's combing separates the fixed part from the moving part }
\label{figure18}
\end{figure}

Recall that the algebraic part of a combed mixed braid is called {\it algebraic mixed braid\/} and it
is an element of $B_{m,n}$. (Recall Figure 6 for an example.) The set $B_{m,n}$ of all algebraic mixed
braids on $m$ fixed strands and $n$ moving strands is closed under the  usual
concatenation and with respect to inverses. Thus, it is a subgroup of $B_{m+n}$. The set  $C_{m,n}$ of
combed mixed braids is a coset of $B_{m,n}$ in $B_{m+n}$ (cf.\ Proposition 1 in \cite{L2}). Thus, for a
fixed manifold $V$, an element in $B_{m,n}$  represents unambiguously an element in $C_{m,n}$, hence an
oriented link in $V$.   The braid group $B_{m,n}$
embedds naturally into the group $B_{m,n+1}$ and we shall denote by $B_{m,\infty} :=
\bigcup_{n=1}^{\infty} B_{m,n}$ the disjoint union of  all braid groups $B_{m,n}.$

\bigbreak   
We would like to restate the Markov equivalence in Theorems 2 and 3 for parted mixed
braids  {\it in terms of their corresponding algebraic mixed braids\/} after
combing. For this we need to understand how exactly the combing is done and how it affects
the parted braid equivalence moves. 

\smallbreak  

Note that, if we regard  a parted mixed braid as an  element of the classical braid group
$B_{m+n},$ then the crossings $\sigma_j$ of the moving part commute with the crossings of
the fixed part, so they are not affected by combing. More precisely,  if $\Sigma_k$
denotes the crossing between the $k$th and the $(k+1)$st strand of the fixed subbraid,
then for all $j=1,2,\ldots,n-1$ and $k=1,2,\ldots,m-1$ we have the relations:
\[
\Sigma_k \sigma_j = \sigma_j \Sigma_k. 
\]  
\noindent Thus the only generating elements
of the moving part that are affected by the combing are the loops $a_i.$ 
 In Lemma 6 below we give formuli for the effect of combing on the $a_i$'s.

\begin{lem}{ \ The crossings $\Sigma_k,$ for $k=1,\ldots,m-1,$ and the loops $a_i,$  for
$i=1,\ldots,m,$ satisfy the following {\rm `combing' relations}:
\[
\begin{array}{llcll}
\bullet \ \  & \Sigma_k {a_k}^{\pm 1}  & = & {a_{k+1}}^{\pm 1} \Sigma_k  &    \\ 

\bullet \ \  & \Sigma_k {a_{k+1}}^{\pm 1}  & = &  {a_{k+1}}^{-1} {a_k}^{\pm 1} a_{k+1}
\Sigma_k  &   \\ 

\bullet \ \  & \Sigma_k {a_i}^{\pm 1}   & = &  {a_i}^{\pm 1}  \Sigma_k  & \mbox{if  \ }
 i \neq k, k+1 \\ 

\bullet \ \  & {\Sigma_k}^{-1} {a_k}^{\pm 1}  & = & a_k {a_{k+1}}^{\pm 1} {a_k}^{-1}
 {\Sigma_k}^{-1}  &    \\ 

\bullet \ \  & {\Sigma_k}^{-1} {a_{k+1}}^{\pm 1}  & = & {a_k}^{\pm 1} 
{\Sigma_k}^{-1}  &   \\ 

\bullet \ \  & {\Sigma_k}^{-1} {a_i}^{\pm 1}   & = &  {a_i}^{\pm 1}  {\Sigma_k}^{-1}  
 & \mbox{if  \ }  i \neq k, k+1.  \\ 
 \end{array}
\]
\noindent Moreover, since $B$ is assumed to be a pure braid for $V= \chi (S^3, \widehat
B),$ it is useful to give the `combing' relations between  the crossings ${\Sigma_k}^2$ 
 and the loops $a_i.$ 
Indeed we have: 
\[
\begin{array}{llcll}
\bullet \ \  & {\Sigma_k}^2 {a_k}^{\pm 1}  & = & {a_{k+1}}^{-1} {a_k}^{\pm 1} a_{k+1} 
{\Sigma_k}^2  &    \\ 

\bullet \ \  & {\Sigma_k}^2 {a_{k+1}}^{\pm 1}  & = &  {a_{k+1}}^{-1} {a_k}^{-1} 
{a_{k+1}}^{\pm 1} a_k a_{k+1} {\Sigma_k}^2  &  \\ 

\bullet \ \  & {\Sigma_k}^2 {a_i}^{\pm 1}   & = &  {a_i}^{\pm 1}  {\Sigma_k}^2  & \mbox{if  \ } 
i \neq k, k+1  \\ 

\bullet \ \  & {\Sigma_k}^{-2} {a_k}^{\pm 1}  & = & a_k a_{k+1}  {a_k}^{\pm 1}
 {a_{k+1}}^{-1} {a_k}^{-1} {\Sigma_k}^{-2}  &    \\ 

\bullet \ \  & {\Sigma_k}^{-2} {a_{k+1}}^{\pm 1}  & = & a_k {a_{k+1}}^{\pm 1} {a_k}^{-1} 
 {\Sigma_k}^{-2}  &   \\ 

\bullet \ \  & {\Sigma_k}^{-2} {a_i}^{\pm 1}   & = &  {a_i}^{\pm 1}  {\Sigma_k}^{-2}  
 & \mbox{if  \ }   i \neq k, k+1.  \\ 

 \end{array}
\]
 }
\end{lem} 

\begin{proof}
  We illustrate in Figure 19 the first three principal relations for
$\Sigma_1.$ For  arbitrary $\Sigma_k$ the proof is obviously analogous. The relations
for the  crossings ${\Sigma_k}^{-1}$ and for ${\Sigma_k}^2$ follow easily from the ones for
$\Sigma_k.$
\end{proof}

\bigbreak 

\begin{figure}[h]
\begin{center}
\includegraphics[width=5in]{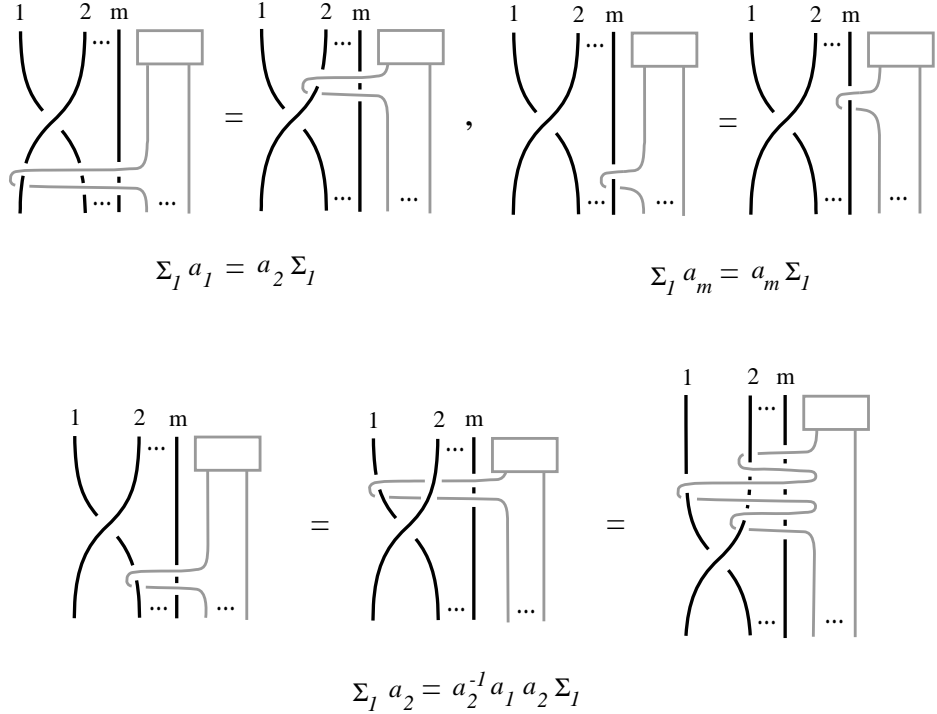}
\end{center}
\caption{ The relations between $\Sigma_1$ and the $a_i$'s }
\label{figure19}
\end{figure}

In $B_{m,\infty}$ we define now the following moves.

\begin{defn}{\rm 
{\it (1) \  Twisted loop conjugation} is defined to be a combed loop conjugation and it has
the algebraic expressions: 
\[
 \beta  \sim  {a_i}^{\mp 1} \beta {\rho_i}^{\pm 1}
\] 
 for  $\beta, a_i, \rho_i \in B_{m,n}$,
where $\rho_i$ is the combing of the loop  $a_i$ through the fixed braid $B$. (Note that the
combing of the loop  ${a_i}^{-1}$ through the fixed braid $B$ is ${\rho_i}^{-1}$.) 

\vspace{.05in}
\noindent {\it (2) \  Algebraic Markov conjugation\/} is Markov conjugation between elements of
$B_{m,\infty}$ and it has the algebraic expression:  
\[ 
\alpha \sim {\sigma_j}^{\pm 1} \alpha {\sigma_j}^{\mp 1}
\] 
 where $\alpha, \sigma_j \in B_{m,n}$. 

\vspace{.05in}
\noindent {\it (3)} \ An {\it algebraic $L$--move\/} is a parted $L$--move between elements of
$B_{m,\infty}$. From Figure 15 one can easily derive the following algebraic expressions for 
algebraic $L_o$--moves and  algebraic $L_u$--moves respectively.  
\[ 
\alpha=\alpha_1\alpha_2 \sim
\sigma_i^{-1}\ldots \sigma_n^{-1} \alpha_1 \sigma_{i-1}^{-1}\ldots
\sigma_{n-1}^{-1}\sigma_n^{\pm 1} \sigma_{n-1} \ldots \sigma_i
 \alpha_2 \sigma_n \ldots \sigma_i 
\]  
 \[ 
\alpha=\alpha_1\alpha_2 \sim
\sigma_i\ldots \sigma_n  \alpha_1 \sigma_{i-1}\ldots
\sigma_{n-1}\sigma_n^{\pm 1}\sigma_{n-1}^{-1}\ldots\sigma_i^{-1}
\widetilde{\alpha_2}\sigma_n^{-1}\ldots\sigma_i^{-1} 
\] 
\noindent  where $\alpha_1, \alpha_2 \in B_{m,n}$.  

\vspace{.05in}
\noindent  {\it (4)} \ An {\it algebraic $M$--move\/} is an $M$--move between elements of
$B_{m,\infty}$  and it has the algebraic expression: 
\[
{\alpha}_1 {\alpha}_2 \sim {\alpha}_1{\sigma_n}^{\pm 1}{\alpha}_2
\]  
where  ${\alpha}_1, {\alpha}_2 \in B_{m,n}.$  
 }
\end{defn}

\begin{lem}{ \ Two parted mixed braids that differ by Markov conjugation by some $\sigma_j,$ 
 resp.\ by an $M$--move, resp.\ by a parted $L$--move, after combing they give rise to algebraic
mixed braids that differ by algebraic Markov conjugation by the $\sigma_j$,
 resp.\ by an algebraic $M$--move, resp.\ by an algebraic $L$--move. 
 }
\end{lem} 

\begin{proof}
 As observed earlier,  the crossings of the moving part commute
with the crossings of the fixed part. Thus, Markov conjugation,   the $M$--moves  and the  
parted $L$--moves all commute with combing. Moreover, the two parted
mixed braids are otherwise identical, so they are both combed in exactly the same manner. 
Therefore, after combing, the combed mixed braids as well as their corresponding  algebraic
mixed braids will just differ by algebraic Markov conjugation,  resp.\ an
algebraic $M$--move,  resp.\ an algebraic $L$--move. 
\end{proof}

We are now in a position to restate Theorem 2  in terms of algebraic mixed braids. 

\begin{thm}[Algebraic Markov Theorem for $S^3 \backslash \widehat B$]  \ Two  oriented
links in  $S^3 \backslash \widehat B$ are isotopic if and only if  any two
corresponding algebraic mixed braid representatives in  $B_{m,\infty}$  differ by a finite
sequence of the following moves:

\vspace{.03in}
\noindent (1) \ Algebraic $M$--moves: \ ${\alpha}_1 {\alpha}_2 \sim
{\alpha}_1{\sigma_n}^{\pm 1}{\alpha}_2, \ \ for \ {\alpha}_1,
{\alpha}_2  \in B_{m,n}$, 

\vspace{.03in}
\noindent (2) \ Algebraic Markov conjugation: \ 
$\alpha \sim {\sigma_j}^{\pm 1} \alpha {\sigma_j}^{\mp 1}$, \ \ for \ $\alpha, \sigma_j \in B_{m,n}$,  

\vspace{.03in}
\noindent (3) \ Twisted loop conjugation: \ 
$\beta  \sim  {a_i}^{\mp 1} \beta {\rho_i}^{\pm 1}$, \ \  for \ $\beta, a_i, \rho_i \in B_{m,n}$, where
$\rho_i$ is the combing of the loop  $a_i$ through $B$, 

\bigbreak
\noindent or, equivalently, by a finite sequence of the following moves: 

\vspace{.03in}
\noindent (1$'$) \  algebraic $L$--moves (see algebraic expressions in Definition 6),

\vspace{.03in}
\noindent (2$'$) \ Twisted loop conjugation.
\end{thm}

\begin{proof}
  By Lemma 7, $M$--moves and Markov conjugation get combed to
algebraic $M$--moves and algebraic Markov conjugation. Thus, by Theorem 2, we only have
to observe that conjugating a parted mixed braid by a loop $a_i$ induces after combing the
twisted conjugation on the  level of the corresponding algebraic braids.   Lemma 6
explains how to do efficiently the combing of the loops ${a_i}^{\pm 1}.$ 
\end{proof}

In order to extend Theorem 4 to mixed braids in $\chi (S^3, \widehat B)$ we need
to understand  how a parted  band move is combed through the surgery braid $B$ and to 
give algebraic expressions for parted band moves between algebraic mixed braids.

\begin{defn}{\rm \ An  {\it algebraic  band move\/} is defined to be a  parted  band move
between elements of $B_{m,\infty}$. See Figure 20 for an abstract example. Setting 
\[
\lambda_{n-1} := \sigma_{n-1} \cdots \sigma_1  \mbox{ \ \ and \ \ } 
t_{k,n} := \sigma_n \cdots \sigma_1 a_k {\sigma_1}^{-1} \cdots {\sigma_n}^{-1},
\]
an algebraic band move has the following algebraic expression: 
\[
\beta_1 \beta_2 \ \sim \ \beta'_1 \, {t_{k,n}}^{p_k} \,
{\sigma_n}^{\pm 1}\, \beta'_2,
\]  
\noindent for $\beta_1, \beta_2 \in B_{m,n},$  where $\beta'_1, \beta'_2  \in
B_{m,n+1}$  are the words $\beta_1, \beta_2$ respectively, with the substitutions: 
\[
\begin{array}{lcl}
  {a_k}^{\pm 1} & \longleftrightarrow & {[({\lambda_{n-1}}^{-1}
{\sigma_n}^{2} \lambda_{n-1}) \, a_k]}^{\pm 1}    \\ 

  {a_i}^{\pm 1} & \longleftrightarrow &   ({\lambda_{n-1}}^{-1} {\sigma_n}^{2} 
 \lambda_{n-1}) \, {a_i}^{\pm 1} \, ({\lambda_{n-1}}^{-1} {\sigma_n}^{2}
\lambda_{n-1})^{-1},   \mbox{ \ \ if \ } i < k   \\ 

 {a_i}^{\pm 1}  & \longleftrightarrow & {a_i}^{\pm 1},  \mbox{ \ \ if \ } i > k.    \\ 
 \end{array}
\]
 Moreover, a  {\it combed algebraic  band move\/} is a parted band move that is the composition of an
algebraic band move with the combing of the parallel strand:
\[ 
\beta_1 \beta_2 \ \sim \ \beta'_1 \, {t_{k,n}}^{p_k} \,
{\sigma_n}^{\pm 1}\, \beta'_2 \, r_k
\]
\noindent where   $r_k$  is the  combing
 of the parted parallel strand  to the $k$th surgery strand through the surgery braid.
 } \end{defn}

\noindent  In Figure 20 note that the isotopy of the
little band in the dotted box is treated as `invisible', that is, as identity in the braid
group. 

\bigbreak 

\begin{figure}[h]
\begin{center}
\includegraphics[width=5.5in]{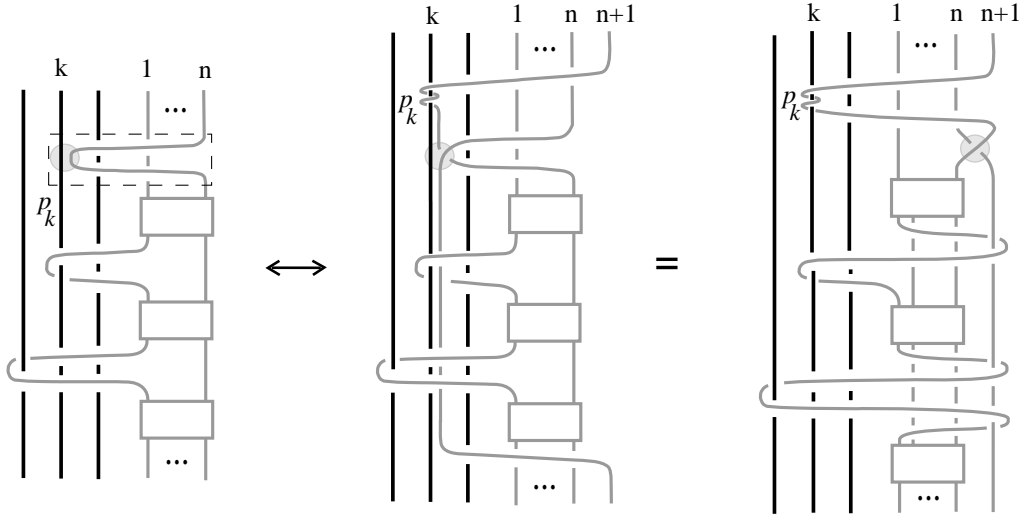}
\end{center}
\caption{  An algebraic band move and its algebraic expression }
\label{figure20}
\end{figure}

\begin{lem}{ \  Performing a parted band move on a parted mixed braid and then combing, the
result is the same as combing the mixed braid and then performing an algebraic band move.
 }
\end{lem} 

\begin{proof}  The parted band move takes place at the top part of the braid, so it resembles an
algebraic band move. Therefore, we just have to consider the behaviour of the parallel strand with
respect to combing. On the other hand, the fact that a band move takes place very close to the surgery
strand ensures that the loops ${a_k}^{\pm 1}$ around the specific surgery strand get combed in the same
way before and after the band move. 

So, when we perform a parted band move on a parted mixed braid we comb away all the
loops ${a_k}^{\pm 1}$ and we leave last the combing of the parallel moving strand. This combing will be
the same in either case of the statement of the Lemma. In Figure 21 we show that using a small braid
isotopy at the bottom of the algebraic part we create an algebraic band move  followed by the combing
through the dotted box $P$ of the parted parallel strand.  Note, finally, that the combing of the
parallel strand leaves clear the fixed braid at the bottom.  Thus the proof is concluded. 
\end{proof}

\bigbreak 

\begin{figure}[h]
\begin{center}
\includegraphics[width=3.7in]{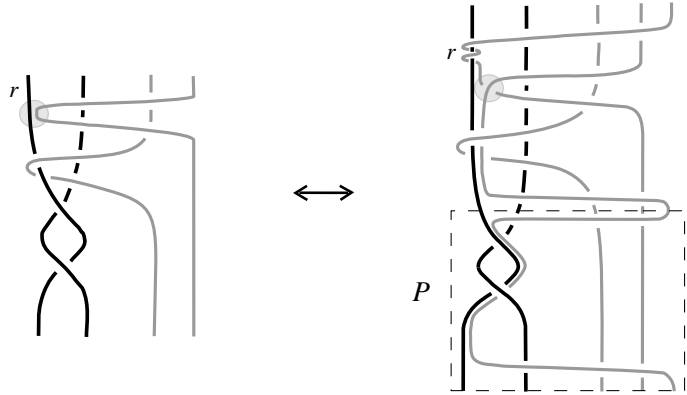}
\end{center}
\caption{  Parted band move =  algebraic band move + combing  }
\label{figure21}
\end{figure}

We are now in the position to state the following result. 

\begin{thm}[Algebraic Markov Theorem for  $\chi (S^3, \widehat B)$]  \ Two  oriented links
in  $\chi (S^3, \widehat B)$ are isotopic if and only if any two
corresponding algebraic mixed braid representatives in  $B_{m,\infty}$  differ by a finite
sequence of the following moves:

\vspace{.03in}
\noindent (1) \ Algebraic $M$--moves: \ ${\alpha}_1 {\alpha}_2 \sim
{\alpha}_1{\sigma_n}^{\pm 1}{\alpha}_2, \ \ for \ {\alpha}_1,
{\alpha}_2  \in B_{m,n}$, 

\vspace{.03in}
\noindent (2) \ Algebraic Markov conjugation: \ 
$\alpha \sim {\sigma_j}^{\pm 1} \alpha {\sigma_j}^{\mp 1}$, \ \ for \ $\alpha, \sigma_j \in B_{m,n}$,  

\vspace{.03in}
\noindent (3) \ Twisted loop conjugation: \ 
$\beta  \sim  {a_i}^{\mp 1} \beta {\rho_i}^{\pm 1}$, \ \  for \ $\beta, a_i, \rho_i \in B_{m,n}$, where
$\rho_i$ is the combing of the loop  $a_i$ through $B$, 

\vspace{.03in}
\noindent (4) \  Combed algebraic band moves: \ For  for every $k=1,\ldots,m$ we have: 
\[ \beta_1 \beta_2 \ \sim \ \beta'_1 \, {t_{k,n}}^{p_k} \,
{\sigma_n}^{\pm 1}\, \beta'_2 \, r_k,
\]
\noindent where  $\beta_1, \beta_2 \in B_{m,n}$ and  $\beta'_1, \beta'_2  \in B_{m,n+1}$ are as in
Definition 7 and where $r_k$  is the  combing
 of the parted parallel strand  to the $k$th surgery strand through $B$,

\bigbreak
\noindent or, equivalently, by a finite sequence of the following moves: 

\vspace{.03in}
\noindent (1$'$) \  algebraic $L$--moves (see algebraic expressions in Definition 6),

\vspace{.03in}
\noindent (2$'$) \ Twisted loop conjugation,

\vspace{.03in}
\noindent (3$'$) \ Combed algebraic band moves.

\end{thm}

\begin{proof} By Theorem 4, we only have to consider the case where a parted band
move takes place and, by Theorem 3, we we only have to consider the behaviour of a parted
band move with respect to combing. This is done in Lemma 8, in the proof of which it is
also explained that  the combing of the parallel strand gives rise to a combed algebraic
band move on the level of  $B_{m,\infty}.$ 
\end{proof}

\begin{rem}{\rm \ We remark that ${t_{k,n}}^{p_k}$ in Definition 7 of an algebraic band
move is just a Markov conjugate of the  loop
${a_k}^{p_k}$ and that these are the appropriate words for defining inductive Markov traces on
quotient algebras of group algebras of $B_{m,n}.$  Note also that the words in the parentheses of the
substitutions of the loops get significantly simplified if we apply a quadratic relation on the
$\sigma_i$'s. Moreover, in Theorem 5 we obtain the best possible control over the band moves of links in
closed $3$--manifolds, and this is very useful for the study of skein modules of closed $3$--manifolds
\cite{P}.
 }
\end{rem}

\section{Special Cases, Examples}

In this section we give the braid equivalences described in Theorems 4 and 5 for specific
examples of knot complements and closed $3$--manifolds. We also discuss the adaptation of the band
move and braid equivalence for the case where the fixed braid $B$ describing the manifold
is not a pure braid, and we study the example  where $B$ is the closure of a
trefoil.

\bigbreak

$\bullet$ Let  $V$ be the {\bf solid torus or the lens space} $\boldsymbol{L(p,1)}$, for a framing $p\in
{\Bbb Z}.$  Then the description of $V$ in $S^3$ is the unknot, and so for $t:=a_1$ and $t_n
:= \sigma_n \ldots \sigma_1 \, t \, {\sigma_1}^{-1} \ldots {\sigma_n}^{-1}$ we have:
\bigbreak 

{\it Two 
oriented links in a {\rm solid torus} are isotopic if and only if 
any two corresponding  mixed braids in  $B_{1,\infty}$  differ by a finite
sequence of the following moves:

\smallbreak
\noindent {\it (1)} \ {\it Algebraic $M$--moves: \ $\alpha \sim \alpha {\sigma_n}^{\pm 1},  \ \ \ \
\alpha \in B_{1,n}$} 

\vspace{.03in}
\noindent  {\it (2)} \ {\it Algebraic Markov conjugation: \ 
$\alpha \sim {\sigma_i}^{\mp 1} \alpha {\sigma_i}^{\pm 1}, \ \ \ \ \alpha, \sigma_i \in B_{1,n}$}  

\vspace{.03in}
\noindent  {\it (3)} \ {\it Loop conjugation: \ 
$\beta \sim {t}^{\mp 1} \beta {t}^{\pm 1}, \ \ \ \ \beta \in B_{1,n}.$ }

\vspace{.05in}
\noindent  Moreover, if the two links lie in $ L(p,1)$ then the
corresponding algebraic mixed braids differ by a finite sequence of the above
moves together with the following:
\smallbreak

\noindent {\it (4)} \ {\it algebraic band moves: } \ For $\beta \in B_{1,n}$ we have: 
\[
\beta \ \sim \ {t_n}^{p}  \, {\sigma_n}^{\pm 1}\, \beta',
\]  
\noindent  where $\beta' \in B_{1,n+1}$  is the word $\beta$ with the substitutions: 
\[
{t}^{\pm 1}  \ \longleftrightarrow  \ {[({\lambda_{n-1}}^{-1}
{\sigma_n}^{2} \lambda_{n-1})  \, t]}^{\pm 1}.  
\]
}
\begin{rem}{\rm \ Constructing all analogues of the 2--variable Jones polynomial in the solid torus via
braids has been completely studied, see  \cite{L1} and references therein. These invariants are related
to the 3rd skein module of the solid torus  \cite{HK}, \cite{T}. Moreover, this last move (4) is used in
\cite{LP} in order to investigate the 3rd skein module of the lens spaces $ L(p,1).$ }
\end{rem}

\noindent $\bullet$  {\bf $V =$  the complement of the $\boldsymbol{m}$--unlink or a connected
sum of
$\boldsymbol{m}$ lens spaces of  type} $\boldsymbol{L(p,1)}$. Then  the fixed braid representing
$V$ is the identity braid, $I_m,$ and so we have:

\vspace{.05in}
\noindent {\it Two oriented links in the {\rm complement of the $m$--unlink} are isotopic if and
only if any two corresponding  mixed braids in  $B_{m,\infty}$  differ by a finite
sequence of the following moves:

\smallbreak
\noindent {\it (1)} \ {\it Algebraic $M$--moves: \ $\alpha \sim \alpha {\sigma_n}^{\pm 1},  \ \ \ \
\alpha  \in B_{m,n}$} 

\vspace{.03in}
\noindent {\it (2)} \ {\it Algebraic Markov conjugation: \ 
$\alpha \sim {\sigma_i}^{\mp 1} \alpha {\sigma_i}^{\pm 1}, \ \ \ \ \alpha, \sigma_i \in B_{m,n}$}  

\vspace{.03in}
\noindent {\it (3)} \ {\it Algebraic loop conjugation: \ 
$ \beta \sim {a_i}^{\mp 1} \beta {a_i}^{\pm 1}, \ \ \ \ \beta \in B_{m,n}, \ i=1,\ldots, m.$}
\smallbreak
(Compare with Theorem 5 in \cite{HL} about braid equivalence in handlebodies.)

\vspace{.05in}
\noindent  Moreover, if the two links lie in the {\rm  connected sum}
$L(p_1,1)\#\cdots \# L(p_m,1),$ where $p_1,\ldots,p_m \in {\Bbb Z},$ then the corresponding
algebraic mixed braids differ  by a finite sequence of the above moves together with the
following:
\smallbreak

\noindent {\it (4)} \ {\it algebraic band moves: } \ For $\beta \in B_{m,n}$  and for
$k=1,\ldots,m$ we have: 
\[
\beta \ \sim \ {t_{k,n}}^{p_k} \, {\sigma_n}^{\pm 1}\,
\beta',
\]  
\noindent  where $\beta' \in B_{m,n+1}$  is the word
$\beta$ with the substitutions: 

\vspace{.03in}
\noindent $\begin{array}{lcl}
{a_k}^{\pm 1} & \longleftrightarrow & {[({\lambda_{n-1}}^{-1}
{\sigma_n}^{2} \lambda_{n-1}) \, a_k]}^{\pm 1}    \\

 {a_i}^{\pm 1} & \longleftrightarrow &   ({\lambda_{n-1}}^{-1} {\sigma_n}^{2} 
 \lambda_{n-1}) \, {a_i}^{\pm 1} \, ({\lambda_{n-1}}^{-1} {\sigma_n}^{2}
\lambda_{n-1})^{-1},   \mbox{ \ \ if \ } i < k   \\ 

{a_i}^{\pm 1}  & \longleftrightarrow & {a_i}^{\pm 1},   \ { \ \ \mbox if} \ i > k.  \\ 
 \end{array}$ }
\bigbreak

\noindent $\bullet$ {\bf  $V =$   the complement of the  Hopf link or a lens space 
$\boldsymbol{L(p,q)}$ obtained by doing surgery along the Hopf link with framings }
 $\boldsymbol{p_1, p_2 \in {\Bbb Z}}$ (obtained from the numerical equation \, $p/q= p_1 + 1/p_2$). The
fixed braid representing $V$ is  ${\Sigma_1}^{2}$ and we have:

\vspace{.03in}
\noindent - {\it Relations for the twisted conjugation:}
\[
\begin{array}{lcl}

 {\Sigma_1}^{2} \cdot {a_1}^{\pm 1} & = & {a_2}^{-1} {a_1}^{\pm 1}a_2 \cdot
{\Sigma_1}^{2},   \\ 

{\Sigma_1}^{2} \cdot {a_2}^{\pm 1} 
& = &  {a_2}^{-1}{a_1}^{-1} {a_2}^{\pm 1} a_1 a_2  \cdot {\Sigma_1}^{2}.   \\ 

 \end{array}
\]
\noindent - {\it Combed algebraic band moves: } \ For  $\beta_1, \beta_2 \in B_{2,n}$ and $k=1,2$ we
have:
\[ 
\beta_1 \beta_2 \ \sim \ \beta'_1 \, {t_{k,n}}^{p_k} \,
{\sigma_n}^{\pm 1}\, \beta'_2 \, r_k
\]
\noindent where  $r_k$  is the  combing through the fixed
braid of the parted moving strand parallel to the $k$th surgery
strand. For $\lambda_n := \sigma_n \cdots \sigma_1$, $r_1, r_2$ are given by the relations:
\[ 
\begin{array}{lcl}
 r_1 & = & \lambda_n {a_2} {\lambda_n}^{-1} \\ 

 r_2 & = & \lambda_n {a_2}^{-1}{a_1}{a_2}{\lambda_n}^{-1}. \\ 
 \end{array}
\]
\noindent  $\beta'_1, \beta'_2 \in B_{2,n+1}$  are the words $\beta_1, \beta_2$ with the following
changes, depending on whether the band move is taking place along the first surgery strand
or along the second. That is, if $k=1$, then $\beta'_1, \beta'_2 $  are obtained from  $\beta_1,
\beta_2$ by doing the substitutions: 
\[ 
\begin{array}{lcl}
 {a_1}^{\pm 1}  & \longleftrightarrow  & {[({\lambda_{n-1}}^{-1}
{\sigma_n}^{2} \lambda_{n-1}) \, a_1]}^{\pm 1}  \\

{a_2}^{\pm 1}  & \longleftrightarrow & 
{a_2}^{\pm 1}.  \\ 
 \end{array}
\]
\noindent  If $k=2$, then $\beta'_1, \beta'_2 $  are obtained from  $\beta_1,
\beta_2$ by doing the substitutions: 
\[
\begin{array}{lcl}

 \ \ \ \ \ \ \ \ \ \ \ {a_1}^{\pm 1} & \longleftrightarrow &   ({\lambda_{n-1}}^{-1}
{\sigma_n}^{2} 
 \lambda_{n-1}) \, {a_1}^{\pm 1} \, ({\lambda_{n-1}}^{-1} {\sigma_n}^{2}
\lambda_{n-1})^{-1}   \\ 

 \ \ \ \ \ \ \ \ \ \ \ {a_2}^{\pm 1} & \longleftrightarrow & {[({\lambda_{n-1}}^{-1}
{\sigma_n}^{2} \lambda_{n-1}) \, a_2]}^{\pm 1}.    \\ 
 \end{array}
\]

\noindent $\bullet$ Let now $V$  be  in general {\bf the complement of a daisy chain on
$\boldsymbol{m}$  rings or a lens space of  type $\boldsymbol{L(p,q)}$ obtained by doing surgery
along the components, with framings} 
$\boldsymbol{p_1, \ldots, p_m \in {\Bbb Z}}$ (which are obtained from the continued fraction expansion of
the rational number $p/q$). The basic manifolds of this series are described in the previous example. A 
fixed braid representing $V$ is 
\[ 
{\Sigma_1}^{2}{\Sigma_3}^{2}\ldots {\Sigma_{2k-1}}^{2}
{\Sigma_2}^{2}{\Sigma_4}^{2}\ldots {\Sigma_{2k-2}}^{2}
 := {\mathcal DC}_{2k},
\]
\noindent  if the daisy chain consists of $m=2k$ rings, and  
\[
{\Sigma_1}^{2}{\Sigma_3}^{2}\ldots {\Sigma_{2k-1}}^{2} 
{\Sigma_2}^{2}{\Sigma_4}^{2}\ldots {\Sigma_{2k}}^{2}
 := {\mathcal DC}_{2k+1},
\]
\noindent  if the daisy chain consists of $m=2k+1$ rings.  It is easy to verify the above braid
words by closing  the
odd-numbered strands by simple arcs that run under the braid and the even-numbered strands 
 by simple arcs that run over the braid.

\bigbreak
\noindent - {\it Relations for the twisted conjugation:} \ We give relations for
$2k$  and $2k+1$ rings by inductive formulas. For $m=2$ the relations for  ${a_1}^{\pm 1}$
and ${a_2}^{\pm 1}$ are given in the previous example.  For
$m=3$  we have the `twisted' relations:

\bigbreak
\noindent $\begin{array}{lcl}
 [{\Sigma_1}^{2} {\Sigma_2}^{2}] \cdot {a_1}^{\pm 1}  & =  &  {a_2}^{-1}
{a_1}^{\pm 1} a_2  \cdot [{\Sigma_1}^{2} {\Sigma_2}^{2}],  \\

[{\Sigma_1}^{2} {\Sigma_2}^{2}] \cdot {a_2}^{\pm 1}  & =  &  (a_1 a_2 a_3)^{-1} {a_2}^{\pm 1} (a_1 a_2
a_3) \cdot [{\Sigma_1}^{2} {\Sigma_2}^{2}],  \\

[{\Sigma_1}^{2} {\Sigma_2}^{2}] \cdot {a_3}^{\pm 1}  & =  & ({a_2}^{-1} {a_1}^{-1} a_2 a_1 a_2
a_3)^{-1} {a_3}^{\pm 1} ({a_2}^{-1} {a_1}^{-1} a_2
a_1 a_2 a_3) \cdot [{\Sigma_1}^{2} {\Sigma_2}^{2}].  \\
 \end{array}$

\bigbreak
\noindent Notice that  the relation for ${a_1}^{\pm 1}$ is the same as for two rings.

\vspace{.05in}
\noindent Suppose now that the twisted conjugation relations are known for $2k-2$  and
$2k-1$ rings and consider $\boldsymbol{2k}$ {\bf rings}. For $1\leq i \leq 2k-3$ the twisted relations
for  ${\mathcal DC}_{2k} \cdot {a_i}^{\pm 1}$ are the same as those for ${\mathcal DC}_{2k-2}
\cdot  {a_i}^{\pm 1}.$  For $i=2k-2,\, 2k-1,\, 2k$ we have the following relations, that
are easy consequences of Lemma 6.

\bigbreak
\noindent $\begin{array}{lcl}
 {\mathcal DC}_{2k} \cdot {a_{2k-2}}^{\pm 1}  & =  & (a_{2k-3} a_{2k-2} {a_{2k}}^{-1} a_{2k-1}
a_{2k})^{-1}  {a_{2k-2}}^{\pm 1} 
(a_{2k-3} a_{2k-2} {a_{2k}}^{-1} a_{2k-1} a_{2k}) \cdot   \\ 

 & & {\mathcal DC}_{2k}, \\

 {\mathcal DC}_{2k} \cdot {a_{2k-1}}^{\pm 1}  & =  & (a_{2k} {a_{2k-2}}^{-1} {a_{2k-3}}^{-1} a_{2k-2} 
a_{2k-3} a_{2k-2} {a_{2k}}^{-1} a_{2k-1} a_{2k})^{-1} \cdot   \\ 

 & & (a_{2k} {a_{2k-2}}^{-1} {a_{2k-3}}^{-1} a_{2k-2} 
a_{2k-3} a_{2k-2} {a_{2k}}^{-1} a_{2k-1} a_{2k}) \cdot {\mathcal DC}_{2k}, \\ 

 {\mathcal DC}_{2k} \cdot {a_{2k}}^{\pm 1}  & =  & (a_{2k-1} a_{2k})^{-1} 
 {a_{2k}}^{\pm 1}  (a_{2k-1} a_{2k}) \cdot {\mathcal DC}_{2k}. \\ 
 \end{array}$

\bigbreak
\noindent Finally, if we have $\boldsymbol{2k+1}$ {\bf rings}, then for $1\leq i \leq 2k-1$
the twisted relations for
${\mathcal DC}_{2k+1} \cdot {a_i}^{\pm 1}$ are the same as those for ${\mathcal DC}_{2k} \cdot
{a_i}^{\pm 1}.$ For $i=2k, 2k+1$ we have the following relations, that are also 
easy consequences of Lemma 6.

\bigbreak
\noindent $\begin{array}{lcl}
 {\mathcal DC}_{2k+1} \cdot {a_{2k}}^{\pm 1}  & =  & 

(a_{2k-1} a_{2k} a_{2k+1})^{-1} {a_{2k}}^{\pm 1}  
(a_{2k-1} a_{2k} a_{2k+1}) \cdot {\mathcal DC}_{2k+1},  \\ 

 {\mathcal DC}_{2k+1} \cdot {a_{2k+1}}^{\pm 1}  & =  & ({a_{2k}}^{-1} {a_{2k-1}}^{-1} a_{2k} a_{2k-1}
a_{2k} a_{2k+1})^{-1} {a_{2k+1}}^{\pm 1} \cdot \\

 & & 
({a_{2k}}^{-1} {a_{2k-1}}^{-1} a_{2k} a_{2k-1} a_{2k} a_{2k+1}) \cdot {\mathcal DC}_{2k+1}.  \\ 
 \end{array}$

\bigbreak
\noindent - {\it Combed algebraic band moves: } \ For  $\beta_1, \beta_2 \in B_{m,n}$  and for
$s=1,\ldots,m$ we have: 
\[ 
\beta_1 \beta_2 \ \sim \ \beta'_1 \, {t_{s,n}}^{p_s} \, {\sigma_n}^{\pm 1}\,
\beta'_2 \, r_s 
\]  
\noindent where $\beta'_1, \beta'_2 \in B_{m,n+1}$  are the words $\beta_1, \beta_2$ with the
substitutions: 

\vspace{.05in}
\noindent $\begin{array}{lcl}

{a_s}^{\pm 1} & \longleftrightarrow & {[({\lambda_{n-1}}^{-1}
{\sigma_n}^{2} \lambda_{n-1}) \, a_s]}^{\pm 1}    \\ 

{a_i}^{\pm 1} & \longleftrightarrow &   ({\lambda_{n-1}}^{-1} {\sigma_n}^{2} 
 \lambda_{n-1}) \, {a_i}^{\pm 1} \, ({\lambda_{n-1}}^{-1} {\sigma_n}^{2}
\lambda_{n-1})^{-1},   \mbox{ \ \ if \ } i < s   \\ 

{a_i}^{\pm 1}  & \longleftrightarrow & {a_i}^{\pm 1},   \mbox{ \ \ if \ } i > s  \\ 
 \end{array}$

\bigbreak
\noindent and where  $r_s \in B_{m,n+1}$  is the  combing through ${\mathcal DC}_{m}$ of the
parted moving strand parallel to the $s$th surgery strand. For  any  index $m$  
of ${\mathcal DC}_{m}$ the combings $r_1, \ldots, r_m$
are given by the following relations: 
\[ 
\begin{array}{lcl}
 r_1 & = & \lambda_n  {a_2} {\lambda_n}^{-1}.   \\ 
 \end{array}
\]
\noindent The combings $r_2, \ldots, r_{m-2}$ are given by the following paired formulas:
\[ 
\begin{array}{lcl}
r_{2k} & = & \lambda_n \, 
({a_{2k}}^{-1} a_{2k-1} a_{2k} {a_{2k+2}}^{-1} {a_{2k+1}} {a_{2k+2}}) \, {\lambda_n}^{-1}  \\ 

r_{2k+1} & = & \lambda_n \, ({a_{2k+1}}^{-1}{a_{2k+2}} {a_{2k}}^{-1}
{a_{2k-1}}^{-1} {a_{2k}} {a_{2k-1}}{a_{2k}} {a_{2k+2}}^{-1} {a_{2k+1}}
{a_{2k+2}}) \,  {\lambda_n}^{-1}.  \\
 \end{array}
\]
\noindent The final combings $r_{m-1}, \, r_m$ depend on whether $m$ is even or odd. For 
 $\boldsymbol{m}$ {\bf even} we have:
\[ 
\begin{array}{lcl}
 r_{m-1} & = & \mbox{as above for odd index}   \\ 
 r_m & = & \lambda_n \, ({a_m}^{-1}{a_{m-1}}{a_m}) \, {\lambda_n}^{-1}.   \\ 
 \end{array}
\]
\noindent For $\boldsymbol{m}$ {\bf odd} we have:
\[ 
\begin{array}{lcl}
 r_{m-1} & = &   \lambda_n \, ( {a_{m-1}}^{-1} a_{m-2} a_{m-1} a_m ) \,
{\lambda_n}^{-1} \\
r_m & = & \lambda_n \,  ({a_m}^{-1} {a_{m-1}}^{-1} {a_{m-2}}^{-1} a_{m-1} {a_{m-2}}{a_{m-1}} a_m) \,
{\lambda_n}^{-1}.  \\ 
 \end{array}
\]

\noindent $\bullet$ {\bf  $V =$  the complement of the Borromean rings or a closed manifold obtained by
doing surgery  along them,  with framings $\boldsymbol{p_1, p_2, p_3 \in {\Bbb Z}}$.} In particular,
with framings
$+1$ we obtain dodecahedral space. The fixed braid representing $V$ is 
\[
{\Sigma_1}^{-1}\Sigma_2 {\Sigma_1}^{-1}\Sigma_2 {\Sigma_1}^{-1}\Sigma_2 := {\mathcal BR} 
\]
\noindent - {\it Relations for the twisted conjugation:}
\[
\begin{array}{lcl}
 {\mathcal BR} \cdot {a_1}^{\pm 1} 
& = &  (a_3 a_1 a_2 {a_1}^{-1} {a_3}^{-1} a_1 {a_2}^{-1}{a_1}^{-1})^{-1} {a_1}^{\pm 1}
(a_3 a_1 a_2 {a_1}^{-1} {a_3}^{-1} a_1 {a_2}^{-1}{a_1}^{-1}) \cdot {\mathcal BR}   \\ 

{\mathcal BR} \cdot {a_2}^{\pm 1} 
& = &  ({a_1}^{-1}{a_3}^{-1}a_1 a_3)^{-1} {a_2}^{\pm 1} ({a_1}^{-1}{a_3}^{-1}a_1 a_3)  \cdot {\mathcal
BR}   \\ 

 {\mathcal BR} \cdot {a_3}^{\pm 1} & = &   (a_1 {a_2}^{-1} {a_1}^{-1} {a_3}^{-1} {a_1}^{-1} a_3 a_1 a_2
{a_1}^{-1} {a_3}^{-1}a_1 a_3)^{-1} {a_3}^{\pm 1} \cdot \\

& & 
 (a_1 {a_2}^{-1} {a_1}^{-1} {a_3}^{-1} {a_1}^{-1} a_3 a_1 a_2 {a_1}^{-1} {a_3}^{-1}a_1 a_3) 
 \cdot {\mathcal BR}   \\ 
 \end{array}
\]
\noindent - {\it Combed algebraic band moves: } \ For  $\beta_1, \beta_2 \in B_{3,n}$  and for
$k=1,2,3$ we have: 
\[ 
\beta_1 \beta_2 \ \sim \ \beta'_1 \, {t_{k,n}}^{p_k} \, {\sigma_n}^{\pm 1}\,
\beta'_2 \, r_k
\]
\noindent where $\beta'_1, \beta'_2 \in B_{3,n+1}$  are the words $\beta_1, \beta_2 $ with the
substitutions: 

\vspace{.05in}
\noindent $\begin{array}{lcl}

{a_k}^{\pm 1} & \longleftrightarrow & {[({\lambda_{n-1}}^{-1}
{\sigma_n}^{2} \lambda_{n-1}) \, a_k]}^{\pm 1}    \\ 

 {a_i}^{\pm 1} & \longleftrightarrow &   ({\lambda_{n-1}}^{-1} {\sigma_n}^{2} 
 \lambda_{n-1}) \, {a_i}^{\pm 1} \, ({\lambda_{n-1}}^{-1} {\sigma_n}^{2}
\lambda_{n-1})^{-1},   \mbox{ \ \ if \ } i < k   \\ 

{a_i}^{\pm 1}  & \longleftrightarrow & {a_i}^{\pm 1},   \ { \ \ \mbox if} \ i > k  \\ 
 \end{array}$

\vspace{.05in}
\noindent and where  $r_k \in B_{3,n+1}$  is the  combing through ${\mathcal BR}$ of the
parted moving strand parallel to the $k$th surgery strand. The combings $r_1, r_2, r_3$
are given by the following relations, which are easy consequences of Lemma 6.
\[ 
\begin{array}{lcl}
 r_1 & = & \lambda_n \, (a_3 a_1 a_2 {a_1}^{-1} {a_3}^{-1} a_1 {a_2}^{-1}
{a_1}^{-1}) \, {\lambda_n}^{-1},  \\ 

 r_2 & = & \lambda_n \, ({a_1}^{-1} {a_3}^{-1} a_1 a_3 ) \,
{\lambda_n}^{-1},  \\ 

 r_1 & = & \lambda_n \, (a_1 {a_2}^{-1} {a_1}^{-1} {a_3}^{-1} {a_1}^{-1}  
a_3 a_1 a_2  {a_1}^{-1} {a_3}^{-1} a_1 a_3) \, {\lambda_n}^{-1}.  \\ 
 \end{array}
\]

{\bf The case of non-pure surgery braids. } \  A closed $3$--manifold may be easier described by
a non-pure surgery braid. As noted in \cite{LR}, Remark 5.11, in this case Theorem 1 
and consequently Theorems 3 and 5 still apply, but now the  band  moves are  more complicated
to express: In this case a band move is modified so that the replacement of the little band
links only with one of the strands of the same surgery component and {\it it runs in parallel to all
other strands of that surgery component}. See Figure 22. In Lemma 9 below we show that such a
band move may be always assumed to have a specific form. We part such a  band move at the top by
pulling all parallel strands to the last positions of the moving part, {\it over} all strands in
between and respecting their order. 

\bigbreak 

\begin{figure}[h]
\begin{center}
\includegraphics[width=3.5in]{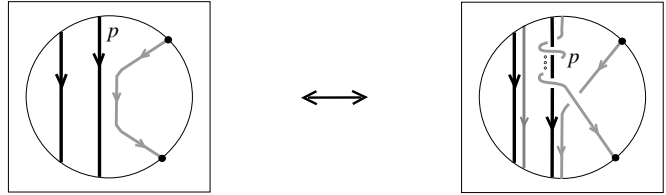}
\end{center}
\caption{  The  band move for a non-pure surgery braid }
\label{figure22}
\end{figure}

\begin{lem}{ \ If the surgery braid $B$ is not a pure braid, then, up to  $L$--equivalence,
a band move may be always assumed to take place at the top part of the braid and on the right of the
rightmost strand of the specific surgery component. 

Moreover, performing a band move at the top part of  a 
mixed braid and then parting, yields the same, up to 
$L$--equivalence and loop conjugation, as performing a parted band move on a parted mixed
braid.
 }
\end{lem}

\begin{proof} Assume that the little band  does not attach to the
rightmost strand of the specific surgery component. Then, as in the proof of Lemma 5,
stretch the little band and its replacement in parallel to the surgery strand and {\it to
the top or to the bottom} of the mixed braid, depending on which direction brings it
 to the  right. See Figure 23.

\bigbreak 

\begin{figure}[h]
\begin{center}
\includegraphics[width=5.3in]{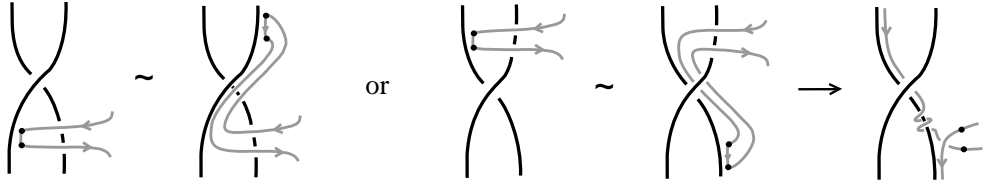}
\end{center}
\caption{  Passing the little band to the right }
\label{figure23}
\end{figure}

If we have landed at the bottom we do appropriate $L$--moves to transfer
the crossing of the band move to the top. Recall Figure 13. If we have reached the
rightmost strand we stop. If not, we continue sliding the little   band and its
replacement till the rightmost strand is reached. Similarly, by loop conjugations we bring the
framing twists to the top of the rightmost strand of the specific component. 

If, moreover, the
two mixed braids that differ by the band move are parted we apply the same strategy as in the
proof of Lemma 5 but adapted here. 
\end{proof}

\noindent $\bullet$ As an application of the above lemma consider $V$ to be {\bf the
complement of the right-handed trefoil or a closed manifold obtained by doing
surgery  along it with framing $\boldsymbol{k \in {\Bbb Z}}$.} In particular,
with framing $-1$ we obtain the dodecahedral space. The fixed braid representing $V$ is ${\Sigma_1}^3$
and we have:

\vspace{.05in}
\noindent - {\it Relations for the twisted conjugation:}
\[
\begin{array}{lcl}
  {\Sigma_1}^3 \cdot {a_1}^{\pm 1}  & = & (a_1 a_2)^{-1} {a_2}^{\pm 1} (a_1 a_2)
\cdot {\Sigma_1}^3     \\ 

  {\Sigma_1}^3 \cdot {a_2}^{\pm 1}  & = & (a_2 a_1 a_2)^{-1} {a_1}^{\pm 1}
(a_2 a_1 a_2) \cdot {\Sigma_1}^3     \\ 
 \end{array}
\]
\noindent - {\it Combed algebraic band moves:} \ For $\beta_1, \beta_2 \in B_{2,n}$ we have:
\[ 
\beta_1 \beta_2 \ \sim \ \beta'_1 \, {\sigma_{n+1}}^{-1} \,  {t_{2,n}}^{k} \, {\sigma_n}^{\pm 1}\,
\sigma_{n+1} \, \beta'_2 \, r
\] 
\noindent where \ $r =   \lambda_n (a_2) \sigma_{n+1} \lambda_n \, ({a_2}^{-1}a_1 a_2)
\, {\sigma_1}^2 {\lambda_n}^{-1} {\sigma_{n+1}}^{-1}
\, ({a_2}^{-1}{a_1}^{-1} a_2 a_1 a_2) \, {\lambda_n}^{-1} \sigma_{n+1}$  

\vspace{.03in}
\noindent is the  combing through the fixed braid
 of the parted moving strands parallel to the two surgery
strands, $t_{2,n} = \lambda_{n} a_k {\lambda_{n}}^{-1}$, and where
$\beta'_1, \beta'_2 \in B_{2,n+2}$  are the words $\beta_1, \beta_2 $ with the substitutions: 
\[
\begin{array}{lcl}
 {a_1}^{\pm 1}  & \longleftrightarrow  & {([{\lambda_{n-1}}^{-1} (\sigma_n
{\sigma_{n+1}}^{2}\sigma_n)
\lambda_{n-1} ] \, a_1 \, [{\lambda_{n-1}}^{-1} (\sigma_n
{\sigma_{n+1}}^{2}{\sigma_n}^{-1})^{-1} \lambda_{n-1}])}^{\pm 1},     \\ 

{a_2}^{\pm 1}  & \longleftrightarrow & {([{\lambda_{n-1}}^{-1} (\sigma_n
{\sigma_{n+1}}^{2}{\sigma_n}^{-1}) \lambda_{n-1}] \, a_2)}^{\pm 1}.    \\ 
 \end{array} 
\]


 \end{document}